\documentclass[11pt, letterpaper, oneside]{amsart}

\usepackage{latexsym, amsmath, amssymb, amsfonts, amscd,bm}
\usepackage{amsthm}
\usepackage{t1enc}
\usepackage[mathscr]{eucal}
\usepackage{indentfirst}
\usepackage{graphicx, pb-diagram}
\usepackage{fancyhdr}
\usepackage{fancybox}
\usepackage{enumerate}
\usepackage[all]{xy}

\usepackage{url}

\theoremstyle{plain}
\newtheorem{thm}{Theorem}[section]

\newtheorem{prop}[thm]{Proposition}

\newtheorem{lemma}[thm]{Lemma}
\newtheorem{cor}[thm]{Corollary}

\theoremstyle{definition}
\newtheorem{defi}[thm]{Definition}

\theoremstyle{remark}
\newtheorem{remark}[thm]{Remark}
\newtheorem{ex}[thm]{Example}

\newcommand{\ZZ}{\ensuremath{\mathbb Z}}

\newcommand{\RR}{\ensuremath{\mathbb R}}
\newcommand{\g}{\ensuremath{\mathfrak{g}}}

\newcommand{\cB}{\mathcal{B}}

\newcommand{\cL}{\mathcal{L}}

\newcommand{\cC}{\mathcal{C}}

\newcommand{\cO}{\mathcal{O}}

\newcommand{\cS}{\mathcal{S}}

\newcommand{\cH}{\mathcal{H}}

\newcommand{\pd}[1]{\frac{\partial}{\partial #1}} 

\newcommand{\ra}{\rangle}
\newcommand{\la}{\langle}
\newcommand{\oo}{[\![}
\newcommand{\cc}{]\!]}
\newcommand{\op}{\textbf{\{}}
\newcommand{\cp}{\textbf{\}}} 

\begin{document}
\title[$L_{\infty}$-algebras and higher analogues]{$L_{\infty}$-algebras and higher analogues of  Dirac structures and Courant algebroids}  
\author{Marco Zambon}
\address{ICMAT(CSIC-UAM-UC3M-UCM)\\
and Departamento de Matem\'aticas,
Universidad Aut\'onoma de Madrid,
Campus de Cantoblanco,
28049 - Madrid, Spain}
\email{marco.zambon@uam.es}


\begin{abstract}
We define    a higher analogue of Dirac structures on  a manifold $M$. Under a regularity assumption, higher Dirac structures can be described by a foliation and a (not necessarily closed, non-unique)   differential form on $M$, and are equivalent to (and simpler to handle than) the Multi-Dirac structures recently introduced in the context of field theory by Vankerschaver, Yoshimura and Marsden.

We associate an $L_{\infty}$-algebra of   observables to every  higher Dirac structure, extending  work of Baez, Hoffnung and Rogers on multisymplectic forms.  Further, applying a recent result of Getzler, we associate an $L_{\infty}$-algebra to any manifold   endowed with a closed differential form $H$, via a  
higher analogue of split Courant algebroid twisted by $H$. 
Finally, we study    
the relations between   the $L_{\infty}$-algebras appearing above.\end{abstract}
\thanks{ 
This work was partially supported by CMUP, by FCT through the programs POCTI, POSI and Ciencia 2007, by grants
PTDC/MAT/098770/2008 and
PTDC/MAT/099880/2008 (Portugal),  
MICINN RYC-2009-04065 and MTM2009-08166-E (Spain).}

\maketitle
\tableofcontents

\section{Introduction}\label{intro}

In the Hamiltonian formalism, many classical mechanical systems are described  by   a manifold, which plays the role of phase space, endowed with a symplectic structure  and a choice of Hamiltonian function.
 However symplectic structures are not suitable to describe all classical systems. Mechanical systems with symmetries are described by Poisson structures -- integrable bivector fields -- and system with constraints are described by closed 2-forms. Systems with   both symmetries and constraints  are described using Dirac structures, introduced by Ted Courant in the early 1990s \cite{Cou}. Recall that, given a manifold $M$,  $TM\oplus T^*M$ is endowed a natural pairing on the fibers and a bracket on its space of sections, called (untwisted) Courant bracket. A Dirac structure     is a maximal isotropic and involutive subbundle of $TM\oplus T^*M$.

Given a Dirac manifold $M$, one defines the notion of Hamiltonian function -- in physical terms, an observable for the system --
and shows that the set of Hamiltonian functions is endowed with a Poisson algebra structure.\\

Higher analogues of symplectic structures are given by  multisymplectic structures  \cite{CIDL}\cite{Hammulti} (called $p$-plectic structures in \cite{BHR}), i.e. closed forms $\omega\in \Omega^{p+1}(M)$  such that the bundle map $\tilde{\omega} \colon TM \to \wedge^p T^*M, X \to \iota_X \omega$ is injective. They are suitable to describe certain physical systems arising from classical field theory, as was realized by Tulczyjew in the late 1960s. They are also suitable to describe systems in which particles are replaced by higher dimensional objects such as  strings \cite{BHR}.

The recent work of Baez, Hoffnung and Rogers \cite{BHR} and then Rogers \cite{RogersL} shows that on a $p$-plectic manifold $M$
the observables -- consisting of certain differential forms  -- have naturally the structure of a  Lie $p$-algebra, by which we mean an $L_{\infty}$-algebra \cite{LadaStasheff} concentrated in degrees $-p+1,\dots,0$.   This extends the fact, mentioned above, that the observables of classical mechanics form a Lie algebra (indeed, a Poisson algebra).\\

The first part of the present paper arose from the geometric observation that, exactly as symplectic structures are special cases of Dirac structures, multisymplectic structures are special cases of higher analogues of Dirac structures. More precisely, for every $p\ge 1$ 
we consider  $$E^p:=TM\oplus \wedge^pT^*M,$$ a vector bundle  endowed with a $\wedge^{p-1}T^*M$-valued pairing and a 
 bracket on its space of sections. We regard $E^p$ as a higher analogue of split Courant algebroids. We also consider isotropic, involutive  subbundles of $E^p$. When the latter are Lagrangian, we refer to them as {higher Dirac structures}.

The following  diagram displays the relations between the geometric structures mentioned so far:
\vspace{-0.2cm}
\begin{center}
\includegraphics[scale=.28]{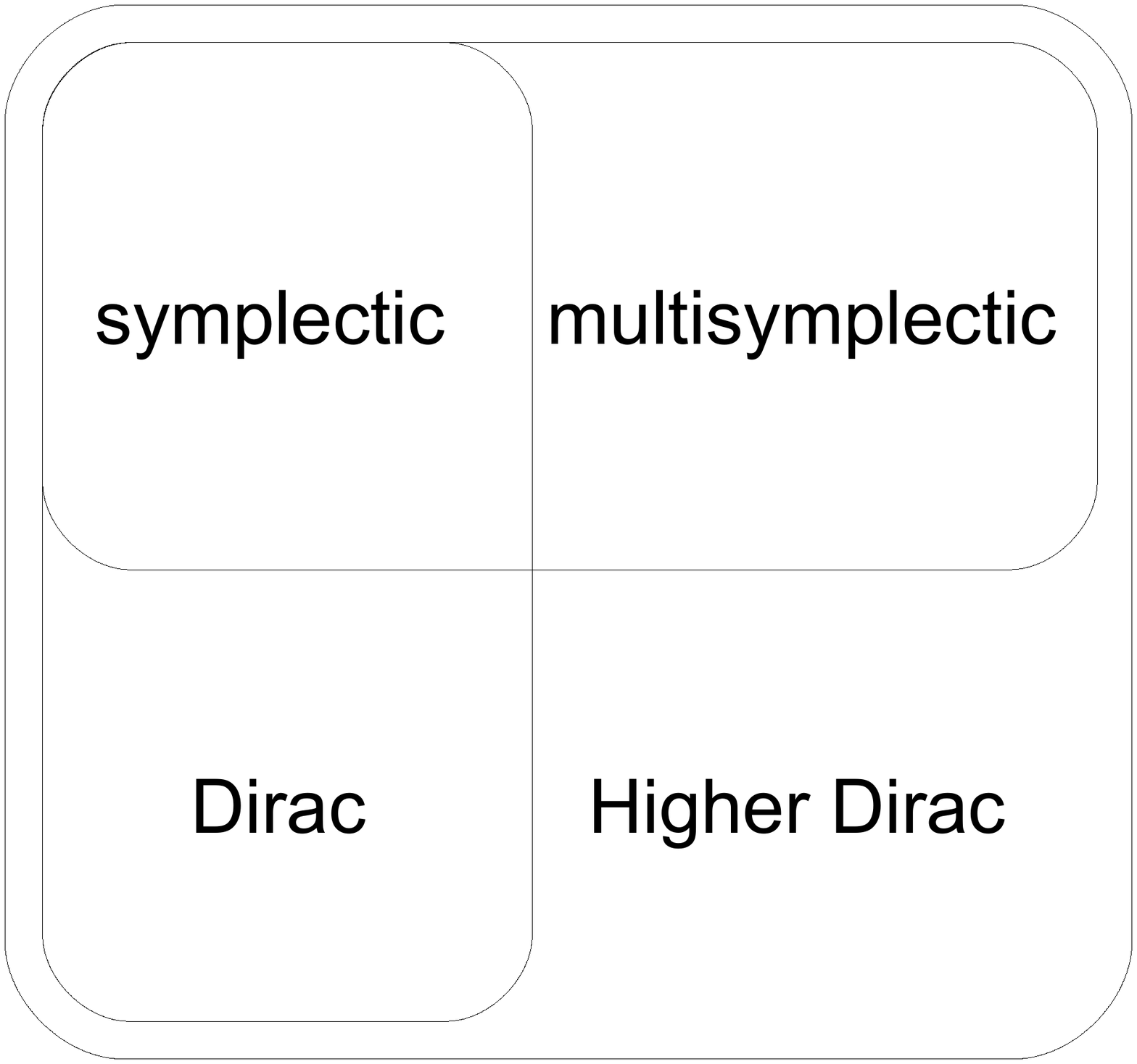}
\end{center}
\vspace{-2.6cm}

In the first part of the paper (\S \ref{haca}-\S \ref{equivmd}) we introduce and study the geometry of isotropic, involutive subbundles of $E^p$. Examples include Dirac structures, closed forms  together with a foliation, and a restrictive class of  multivector fields. 
The main results are
\begin{itemize}
\item Thm. \ref{integ}: a description of all regular higher Dirac structures   in terms of familiar geometric data: a (not necessarily closed) differential form and a foliation. 
\item Thm \ref{eqint}:  higher Dirac structures are \emph{equivalent} to    Multi-Dirac structures, at least in the regular case\footnote{Regularity is a technical assumption and is probably not necessary. The physically most relevant examples of Multi-Dirac structures   are regular \cite{MultiDirac}.}.\end{itemize}
Recall that Multi-Dirac structures were recently introduced by  Vankerschaver, Yoshimura and Marsden \cite{MultiDirac}. They are the geometric structures that allow to describe the implicit Euler-Lagrange equations (equations of motion) of a large class of  field theories, which include the treatment of non-holonomic constraints. 
By the above equivalence, higher Dirac structures thus acquire a field-theoretic motivation. Further, since higher Dirac structures are simpler  to handle than Multi-Dirac structures (which contain some redundancy in their definition), we expect 
our work to be useful in the context of field theory too. \\

The second part of the paper is concerned with the algebraic structure on  the observables, which turns out to be an $L_{\infty}$-algebra. Further, we investigate an $L_{\infty}$-algebra that can be associated to a manifold without any geometric structure on it, except for a (possibly vanishing) closed differential form defining a   
 twist. Recall that a closed 2-form on a manifold $M$ (a 2-cocycle for the Lie algebroid $TM$) can be used to obtain a Lie algebroid structure on $E^0=TM\times \RR$  \cite[\S 1.1]{MariusPre}, so the sections of the latter form a Lie algebra.
 Recall also that Roytenberg and Weinstein \cite{rw} associated  a  Lie 2-algebra to every Courant algebroid (in particular to $E^1=TM\oplus T^*M$ with Courant bracket twisted by a closed 3-form). Recently Getzler  \cite{GetzlerHigherDer} gave an algebraic construction which  extends Roytenberg and Weinstein's proof. Applying Getzler's result in a straightforward way one can extend the above results to all  $E^p$'s.

Our main results in the second part of the paper (\S \ref{Linfty}-\S \ref{per}) are:
\begin{itemize}
\item  Thm. \ref{Liep}: the observables associated to an isotropic, involutive  subbundle of $E^p$ form a Lie $p$-algebra. 
\item Prop. \ref{ord} and Prop. \ref{ordH}: to  $E^p=TM\oplus \wedge^p  T^*M$ and to a closed $p+2$- form $H$ on $M$, one can associate a Lie $p+1$-algebra extending the $H$-twisted Courant bracket. 
\item  
Thm. \ref{mor01}: there is a morphism (with one dimensional kernel) from the Lie algebra associated to $E^0$ and a closed $2$-form   into the  Lie 2-algebra associated  to the Courant algebroid $E^1=TM\oplus T^*M$ with the untwisted Courant bracket.
\end{itemize}
Rogers \cite{RogersCou} observed
that there is an injective  morphism -- which can be interpreted as a  prequantization map -- from the Lie 2-algebra of observables on a $2$-pletic manifold $(M,\omega)$ into the Lie 2-algebra associated   to the Courant algebroid $E^1=TM\oplus T^*M$ endowed with the $\omega$-twisted Courant bracket. We conclude the paper with an attempt to put this into context.\\

\noindent \textbf{Acknowledgments}
I thank  Klaus Bering, Ya\"el Fr\'egier, David Iglesias, Camille Laurent, Jo\~ao Martins, Claude Roger, Chris Rogers, Florian Sch\"atz, Pavol {\v{S}}evera and Joris Vankerschaver for helpful discussions, and Jim Stasheff for comments on this note. The first part of Prop. \ref{referee} on integration is due to a referee, whom I hereby thank. I am grateful to a further referee for numerous comments that improved the presentation.
Further I thank Juan Carlos Marrero and Edith Padr\'on for pointing out to me the reference \cite{Hagi}, and Chris Rogers for pointing out \cite{MultiDirac}.

\section{Higher analogues of split Courant algebroids}\label{haca}

Let $M$ be a manifold and $p\ge 0$ an integer.
Consider  the vector bundle $$E^p:=TM\oplus \wedge^pT^*M,$$
endowed with the 
symmetric pairing on its fibers
\begin{equation*}
\la \cdot,\cdot \ra \colon E^p \times E^p \to \wedge^{p-1}T^*M,
\end{equation*}
given by
\begin{equation}\label{pairing}
\langle X+\alpha , Y+\beta \rangle =  \iota_{X}\beta + \iota_{Y}\alpha .
\end{equation}
Endow the space of sections of $E^p$ with the \emph{Dorfman bracket}
\begin{equation}\label{dorf}
\oo X+\alpha , Y+\beta\cc= [X,Y] + \cL_{X} \beta - \iota_{Y} d\alpha .
\end{equation}

The Dorfman bracket satisfies the Jacobi identity and Leibniz rules
\begin{align}
\label{Jacobi}
\oo e_1,\oo e_2,e_3\cc\cc&=\oo \oo e_1,e_2\cc ,e_3\cc+\oo e_2,\oo e_1,e_3\cc\cc\\
\label{Leib2}
\oo e_1,fe_2 \cc&=f\oo e_1,e_2 \cc+ (pr_{TM}(e_1)f) e_2\\
\label{Leib1}
\oo fe_1,e_2 \cc&=f\oo e_1,e_2\cc- (pr_{TM}(e_2)f) e_1+df\wedge \langle e_1,e_2 \rangle
\end{align}
where $e_i\in\Gamma(E^p)$, $f\in C^{\infty}(M)$, and $pr_{TM} \colon E^p \to TM$ is the projection onto the first factor.
%

 
The decomposition of the Dorfman bracket into its anti-symmetric and symmetric parts is
\begin{equation}\label{dorcou}
\oo e_1,e_2\cc =\oo e_1,e_2\cc_{Cou}+\frac{1}{2}d \la e_1,e_2 \ra,
\end{equation}
where $$\oo X+\alpha , Y+\beta\cc_{Cou}:= [X,Y] + \cL_{X} \beta -\cL_{Y} \alpha -\frac{1}{2}d(\iota_X \beta- \iota_{Y} \alpha)$$  is known as \emph{Courant bracket}.  \\

\begin{remark}
The Dorfman bracket on $E^p$ was already considered by Hawigara \cite[\S 3.2]{Hagi}, Hitchin \cite{Hi} and Gualtieri \cite[\S 3.8]{Gu}\cite[\S 2.1]{Gu2}.
$(E^p,\la \cdot,\cdot \ra , \oo \cdot,\cdot \cc)$ is an example of \emph{weak Courant-Dorfman algebra} as introduced by Ekstrand and Zabzine in \cite[Appendix]{ekzab}.
When $p=1$ we recover  an instance of split  \emph{Courant algebroid} \cite{LWX}. The Courant bracket has been extended to the setting of multivector fields in \cite[\S 4]{MultiDirac}.
\end{remark}

In \cite{Hi,Gu,Gu2} it is remarked that closed $p+1$-forms $B$ on $M$ provide symmetries of the Dorfman bracket (and of the pairing), 
 by the gauge
transformation 
$e^B \colon X+\alpha \mapsto X+\alpha+\iota_XB$.
Further  
the Dorfman bracket may be twisted by a closed  $p+2$-form $H$, just by adding a term 
$\iota_Y\iota_X H$ to the r.h.s. of eq. \eqref{dorf}. We refer to the resulting bracket as \emph{$H$-twisted Dorfman bracket} (this notion will not be used until \S \ref{liM}), and we use  the term Dorfman bracket to refer to the untwisted one given by eq. \eqref{dorf}.

\section{Higher analogues of Dirac structures}\label{hads}

In this section we introduce a geometric structure that extends the notion of Dirac structure and multisymplectic form.  It is given by a subbundle of $E^p$,
which we require to be involutive and  isotropic, since this is needed to associate to it an $L_{\infty}$-algebra of observables  in \S \ref{obs}. Further we consider  subbundles which are  Lagrangian  (that is, maximal isotropic) and study in detail their geometry.

\begin{defi}\label{hd} Let $p\ge 1$. Let $L$ be  a subbundle of $E^p=TM\oplus \wedge^pT^*M$.
\begin{itemize}
\item 
$L$ is \emph{isotropic} 
if for all  sections $X_i+\alpha_i$: 
\begin{equation}\label{isot}
\la X_1+\alpha_1,X_2+\alpha_2\ra=0.
\end{equation}
  $L$ is \emph{involutive} if for all  sections $X_i+\alpha_i$: $$\oo  X_1+\alpha_1,X_2+\alpha_2 \cc \in \Gamma(L),$$
 where $\oo  \cdot,\cdot \cc$ denotes the Dorfman bracket \eqref{dorf}.

\item 
$L$ is \emph{Lagrangian} if  $$L=L^{\perp}:=\{e \in E^p : \la e,L\ra =0  \}.$$

\noindent
(In this case we also refer to $L$  as a \emph{almost Dirac structure of order $p$}.)

\noindent $L$   a \emph{Dirac structure of order $p$} or \emph{higher Dirac structure} if it is Lagrangian and involutive.

\item $L$ is  \emph{regular} if $pr_{TM} (L)$ has constant rank along $M$. 
\end{itemize}


\end{defi}

\subsection{Involutive isotropic subbundles}

In this subsection we make some simple considerations on involutive isotropic subbundles and present some examples.

The involutive, Lagrangian subbundles of $E^1$ are the  \emph{Dirac structures} introduced by Courant \cite{Cou}.

When $p=dim(M)$, isotropic subbundles are forced to lie inside $TM\oplus\{0\}$ or $\{0\}\oplus \wedge^pT^*M$, hence they are uninteresting. 

Now, for arbitrary $p$, we look at involutive, isotropic subbundles that project isomorphically onto the first or second summand of $E^p$.

\begin{prop}\label{closedform} Let $p\ge1$. Let   $\omega$ be a closed $p+1$-form on $M$. Then 
$$graph(\omega):=\{X-\iota_X \omega: X\in TM\}$$ 
is an isotropic involutive subbundle of $E^p$.
All isotropic involutive subbundles  $L\subset E^p$ that project isomorphically onto $TM$ under $pr_{TM} \colon E^p \to TM$ are of the above form. 
\end{prop}
\begin{proof}
The subbundle $graph(\omega)$ is isotropic because $\la X-\iota_X \omega, Y-\iota_Y \omega \ra = -\iota_X\iota_Y \omega-\iota_Y\iota_X \omega=0$.
To see that $L$ is involutive,  use the fact that since $\omega$ is closed $d(\iota_X \omega)=\cL_X\omega$ and compute
$$\oo X-\iota_X \omega, Y-\iota_Y\omega \cc=
[X,Y]-\cL_X(\iota_Y\omega)+\iota_Y(\cL_X\omega)=
[X,Y]- \iota_{[X,Y]}\omega.$$

Let $L\subset E^p$ be a subbundle that  projects isomorphically onto $TM$, i.e. $L=\{X+B(X): X\in TM\}$ for some $B\colon TM \to \wedge^pT^*M$. If $L$ is isotropic then
the map $$TM\otimes TM \to  \wedge^{p-1}T^*M,\;\;\; X\otimes Y \mapsto \iota_X(B(Y))$$ is skew in $X$ and $Y$, so
$B(X)=-\iota_{X}\omega$ defines a unique $p+1$-form $\omega$, which satisfies $graph(\omega)=L$. If $L$ is involutive  then the above computation shows that $\omega$ is a closed form.
\end{proof} 
 
The following   generalization of  Prop. \ref{closedform} is proven exactly as in  the last paragraph of the proof of Prop. \ref{integ}. It provides a wide class of regular isotropic, involutive subbundles.
\begin{cor}\label{oS}
Fix $p\ge 1$. Let $\omega\in \Omega^{p+1}(M)$ be a $p+1$-form  and $S$ an integrable distribution on $M$, such that $d\omega|_{\wedge^3 S \otimes \wedge^{p-1}TM}=0.$ Then
$$L:=\{X-\iota_X \omega +\alpha: X\in S, \alpha\in \wedge^p S^{\circ}\}$$
is an isotropic, involutive subbundle of $E^p$.
\end{cor}
  \begin{prop}\label{NPiso} Let $1 \le p\le dim(M)-1$.  Let $\pi\in \Gamma(\wedge^{p+1} TM)$ be either a Poisson bivector field,  a $dim(M)$-multivector field   or $\pi=0$.
  Then 
$$graph(\pi):=\{\iota_{\alpha} \pi + \alpha: \alpha \in \wedge^{p}T^*M\}$$
is an isotropic involutive subbundle of $E^p$.

All isotropic involutive subbundles  $L\subset E^p$ that project isomorphically onto $\wedge^pT^*M$ under $pr_{\wedge^pT^*M} \colon E^p \to \wedge^pT^*M$ are of the above form. 
\end{prop}
\begin{proof} We write $n:=p+1$, so $\pi$ is an $n$-vector field.
Clearly $graph(\pi)$ is isotropic in the cases $\pi=0$ and $n=2$. For the case $n=dim(M)$ fix a point $x\in M$. We may assume that at $x$ we have  $\pi=\pd{x_1}\wedge \dots \wedge \pd{x_n}$ where $\{x_i\}_{i\le n}$ is a coordinate system on $M$. For each $i$ denote  $dx_i^C:=dx_1\wedge \dots \widehat{dx_i}  \dots \wedge dx_{n}$. 
For $i\le j$ at the point $x$ we have
\begin{align}\label{topm}
&\la \iota_{dx_i^C}\pi+ {dx_i^C},\iota_{dx_j^C}\pi+{dx_j^C} \ra\\
=&((-1)^{(n-i)+(i-1)}+(-1)^{(n-j)+(j-2)})dx_1\wedge \dots \widehat{dx_i} \dots \widehat{dx_j}   \dots \wedge dx_{n}=0  \nonumber
\end{align}
showing that $graph(\pi)$ is isotropic.

It is known that $graph(\pi)$ is involutive if{f} $\pi$ is a Nambu-Poisson multivector field (see \cite[\S 4.2]{Hagi}). For $n=2$ the Nambu-Poisson multivector fields are exactly Poisson bivector field, and for $n=dim(M)$ all $n$-multivector fields are Nambu-Poisson. This concludes the first part of the proof.

Conversely, assume that $L\subset E^{n-1}$ is an isotropic subbundle that projects isomorphically onto $\wedge^{n-1}T^*M$, i.e. 
$L=\{A {\alpha} + \alpha: \alpha \in \wedge^{n-1}T^*M\}$ for some map $A \colon \wedge^{n-1}T^*M \to TM$. 

Assume that  $A$ is not identically zero, and that $n\neq 2, dim(M)$.
In this case we obtain a contradiction to the isotropicity of $L$, as follows.
There is a point $x\in M$ with $A_x\neq 0$. Near $x$ choose coordinates  $x_1,\dots,x_{dim(M)}$ (notice that $dim(M)\ge n+1$).
Without loss of generality at $x$ we might assume that   $A (dx_1\wedge \dots\wedge dx_{n-1})$ does not vanish.
It does not lie in the span of $\pd{x_1}, \dots,\pd{x_{n-1}}$ since we assume that $L$ is isotropic, so by modifying the coordinates  $x_{n},\dots,x_{dim(M)}$
we may assume that $A (dx_1\wedge \dots\wedge dx_{n-1})=\pd{x_{n}}$.
Then $$\big\la A_x(dx_1\wedge \dots \wedge dx_{n-1})+ dx_1\wedge \dots \wedge dx_{n-1}\;,\;
A_x(dx_3\wedge \dots \wedge dx_{n+1})+dx_3\wedge \dots \wedge dx_{n+1} \big\ra \neq 0.$$
Indeed the contraction of $A_x(dx_1\wedge \dots \wedge dx_{n-1})=\pd{x_{n}}$
with $dx_3\wedge \dots \wedge dx_{n+1}$  contains the summand $(-1)^{n-3}\cdot
dx_3\wedge \dots \wedge dx_{n-1}\wedge dx_{n+1}$,
whereas the contraction of any vector of $T_xM$ with $dx_1\wedge \dots \wedge dx_{n-1}$ can not contain $dx_{n+1}$. Hence we obtain a contradiction to the isotropicity.

If $A\equiv0$ then clearly $L$ is isotropic.
In the case $n=2$, it is known that $L$ is isotropic if{f} it is the graph of a bivector field $\pi$.
Now consider the case $n=dim(M)$.
For any $i$, let  $X_i+{dx_i^C}\in L$. The isotropicity condition implies that $X_i=\lambda_i\pd{x_i}$ for some $\lambda_i\in \RR$, and a computation similar to  \eqref{topm} implies $\lambda_i=(-1)^{n-i}\lambda_n$ for all $i$, so that
$L=graph(\pi)$ for   $\pi=\lambda_n\pd{x_1}\wedge \dots \wedge \pd{x_n}$.

Hence we have shown that $L$ is isotropic if{f} $L$ is the graph of an $n$-vector field where $\pi=0$,
$n=2$ or $n=dim(M)$. As seen earlier, if $graph(\pi)$ is involutive then, in the case $n=2$, $\pi$ has to be a Poisson bivector field.
\end{proof}


We present a class of isotropic involutive subbundles   which are not necessarily regular:
\begin{cor}\label{0n} 
Let $\Omega$ be an top degree form on $M$, 
 and $f\in C^{\infty}(M)$ such that $\Omega_x\neq 0$ at points of $\{x\in M: f(x)=0\}$. Then
$$L:=\{fX-\iota_X\Omega: X\in TM\}$$
is an involutive isotropic subbundle of  $E^{dim(M)-1}$.
 \end{cor}
\begin{proof} Let $x\in M$. If $f(x)\neq 0$, then nearby $L$ is the graph of $\frac{1}{f}\Omega$,
 which being a top-form is closed. Hence, near $x$, $L$ defines an isotropic involutive subbundle   by Prop. \ref{closedform}.
Now suppose that  $f(x)=0$. Then $L_x$ is just $0+\wedge^{dim(M)-1}T^*_xM$,
so nearby $L$ is the graph of a top multivector field, and by Prop. \ref{NPiso} it is 
an isotropic involutive subbundle.
\end{proof}

Notice that the isotropic subbundles described in Prop. \ref{closedform}, Prop. \ref{NPiso}, Cor. \ref{0n}  are all Lagrangian (use Lemma \ref{easychar} below).\\



We end this subsection relating  
involutive isotropic subbundles with Lie algebroids and Lie groupoids.

\begin{prop}\label{algoid}
Let $L\subset E^p$ be an involutive isotropic subbundle. Then $(L, \oo \cdot,\cdot \cc,pr_{TM})$ is a Lie algebroid \cite{CW},
where $pr_{TM} \colon E^p\to TM$ is
 the projection onto the first factor.
\end{prop}
\begin{proof}
The restriction of the Dorfman bracket to $\Gamma(L)$ is skew-symmetric because of eq. \eqref{dorcou},
and as seen in eq. \eqref{Jacobi} the Dorfman bracket satisfies the Jacobi identity. The Leibniz rule holds because of  eq. \eqref{Leib2}.
\end{proof}

Recall that (integrable) Dirac structures give rise to presymplectic groupoids in the sense of \cite{BCWZ} and, restricting to the non-degenerate case, that  
Poisson structures give rise to symplectic groupoids. We generalize this:

\begin{prop}\label{referee} Suppose that the 
Lie algebroid $L$ of Prop. \ref{algoid} integrates to a source simply connected Lie groupoid $\Gamma$. Then $\Gamma$ is canonically endowed with a  multiplicative closed $p+1$-form $\Omega$.

Further, if $L$ is the graph of a multivector field as in Prop. \ref{NPiso} or the graph of a 
multisymplectic form (see \S \ref{intro}), then $\Omega$ is a multisymplectic form.
\end{prop}
\begin{proof}
The first statement follows immediately from recent results of  Arias Abad-Crainic, applying \cite[Thm. 6.1]{Abad:2009zr}
to the vector bundle map $\tau \colon L \to \wedge^p T^*M$ given by the projection onto the second factor, which satisfies the assumptions of the theorem since $L$ isotropic and because the Lie algebroid bracket on $L$ is the restriction of the Dorfman bracket. 
Concretely, for all $x\in M$ and $e\in L_x$, $X_1,\dots,X_p \in T_xM$, the multiplicative form $\Omega$ is determined by the equation 
\begin{equation}\label{henrmulti}
\Omega(e,X_1,\dots,X_p)= \langle pr_{\wedge^{p}T^*M}(e), X_1\wedge \dots \wedge X_p \rangle.
\end{equation}
Here on the l.h.s. we identify the Lie algebroid $L$ with $ker(s_*)|_M$, where $s \colon \Gamma \to M$ is the source map.

Now assume that $L$ is the graph of a multivector field $\pi$ as in  Prop. \ref{oS}.
First, given a non-zero $e\in L$, it follows that $pr_{\wedge^{p}T^*M}(e)$ is also non-zero, so it pairs non-trivially with some $X_1\wedge \dots \wedge X_p \in \wedge^{p}TM$. Second, given a non-zero $X_1\in TM$, extend it to a non-zero element $X_1\wedge \dots \wedge X_p
\in \wedge^{p}TM$, and choose $\alpha \in \wedge^{p}T^*M$  so that their pairing is non-trivial. Let $e:=\iota_{  \alpha} \pi+\alpha$. 
Then the expression \eqref{henrmulti} is non-zero. Since $T\Gamma|_M=TM\oplus ker(s_*)|_M$ and $\Omega|_{\wedge^{p+1}TM}=0$, this
 shows that $\Omega$ is multisymplectic at points of $M$. To make the same conclusion at every $g\in \Gamma$, use \cite[eq. (3.4)]{BCWZ}
 that the multiplicativity of $\Omega$ implies \begin{equation*}
\Omega_g((R_g)_* e,w_1,\dots,w_p)=\Omega_{x}(e,t_*(w_1),\dots,t_*(w_p))
\end{equation*}
for all $e\in ker(s_*)|_{x}$ and $w_i\in T_g \Gamma$. Here $t\colon \Gamma \to M$ is the target map and $x:=t(g)\in M$. 

Last, assume that $L$ is the graph of a multisymplectic form $\omega$ on $M$.
Given a non-zero $e\in L$, say $e=X-\iota_X \omega$, we have by eq. \eqref{henrmulti} that $\iota_e \Omega|_{\wedge^{p}TM}=-\iota_X \omega\neq 0$. Given a non-zero $X_1\in TM$, there is $X\wedge X_2\wedge  \dots \wedge X_p
\in \wedge^{p}TM$ with which $\iota_{X_1}\omega$ pairs non-trivially. Let $e:=X-\iota_X \omega$. Then the expression \eqref{henrmulti} is non-zero. This shows that $\Omega$ is multisymplectic at points of $M$, and by the  argument above on the whole of $\Gamma$.
\end{proof}

\subsection{Higher Dirac structures}\label{subsec:lag}

In this subsection we  characterize Lagrangian subbundles $L\subset E^p$ (i.e. almost Dirac structures of order $p$) and their involutivity.

We start characterizing Lagrangian subbundles at the linear  algebra level. Recall first what happens in the case  $p=1$. Let $T$ be a vector space. Any $L\subset T\oplus T^*$ such that $L=L^{\perp}$ is determined exactly by the subspace $S:=pr_T(L)$ and a skew-symmetric bilinear form on it \cite{BR}.
Further $dim(S)$ can assume any value between $0$ and $dim(T)$. For $p\ge 2$
 the description is more involved, however
it remains true that    every Lagrangian subspace of $T\oplus \wedge^p T^*$ can be described by means of a subspace $S\subset T$ (satisfying a dimensional constraint) and a (non-unique) $p+1$-form on $T$.

\begin{prop}\label{linalg}
Fix a vector space $T$ and an  integer $p\ge 1$. There is a bijection between
\begin{itemize}
\item Lagrangian subspaces $L\subset T\oplus \wedge^p T^*$ 
\item \text{ pairs }
$$\begin{cases}
 S\subset T &\text{ such that either } dim(S)\le (dim(T)-p) \text{ or }S=T,\\
 \Omega \in \wedge^2 S^*\otimes \wedge^{p-1}T^* &\text{ such that } 
\Omega \text{ is the restriction of an element of } \wedge^{p+1}T^*.
 \end{cases}$$
\end{itemize}


The correspondence is given by
\begin{align*}
L &\mapsto \begin{cases} S:=pr_T(L)\\
\Omega \text{ given by } \iota_X \Omega
=\alpha|_{S\otimes\bigotimes^{p-1}T} \text{ for all } X+\alpha \in L
\end{cases}\\
(S,\Omega)&\mapsto L:=\{X+\alpha: X\in S, \alpha|_{S\otimes\bigotimes^{p-1}T}=\iota_X \Omega
\}.
\end{align*}
\end{prop}

Here we regard $\wedge^n T^*$ as the subspace of $ \bigotimes^n T^*:=T^*\otimes\dots\otimes T^*$  consisting of elements invariant under the odd representation of the permutation group in $n$ elements. Loosely speaking, the restriction on $dim(S)$ arises as follows: when it is not satisfied $\wedge^p S^{\circ}=\{0\}$ and $S\neq T$,
and one can enlarge 
$L$ to an isotropic $L'\subset  T\oplus \wedge^p T^*$ such that $pr_T(L')$ is strictly larger than $S$. 
The proof of Prop. \ref{linalg} is presented 
in Appendix \ref{applag}.


 
An immediate corollary of    Prop. \ref{linalg}, which we present without proof,  is:
\begin{cor}\label{norom} Fix a vector space $T$ and an  integer $p\ge 1$. For any Lagrangian subspace $L\subset T\oplus \wedge^p T^*$ let $(S,\Omega)$ be the corresponding pair as in Prop. \ref{linalg}, and 
 $\omega\in \wedge^{p+1}T^*$ an arbitrary extension of $\Omega$. Then $L$ can be described in terms of $S$ and $\omega$ as
$$L=\{X+\iota_X \omega +\alpha: X\in S, \alpha\in \wedge^p S^{\circ}\}.$$
\end{cor}

As an immediate consequence of Lemma \ref{easychar}, we obtain the following dimensional constraints on the singular distribution induced by a Lagrangian subbundle:
\begin{cor}\label{restr} Let $L\subset E^p$ be a Lagrangian  subbundle. Denote $S:=pr_{TM}(L)$. Then
\begin{itemize}
\item[a)] $dim(S_x)\in\{0,1,\dots,dim(M)-p,dim(M)\}$ for all $x\in M$
\item[b)] $dim(L_x)=dim(S_x)+ 
{dim(M)-dim(S_x) \choose p}$ is constant for all $x\in M$.
\end{itemize}
 \end{cor}

When $p=1$, so that $L$ is a maximal isotropic subbundle of $TM\oplus T^*M$,
the dimensional constraints of Cor. \ref{restr} do not pose any restriction of $dim(S_x)$. (It is known, however, that 
$dim(S_x) \;mod\; 2$ must be constant on $M$.) When $p\ge 2$,    Lagrangian  subbundles of $E^p$ are quite rigid.
 
\begin{ex}
Let $p=dim(M)-1$, and let $L$ be  a Lagrangian  subbundle of $E^p$. 
Cor. \ref{restr} a) 
implies that at every point $dim(S_x)$ is either  $0$, $1$ or $dim(M)$.
Assume that $p\ge 2$.  
By Cor. \ref{restr} b),  if $rk(S)=1$ at one point then $rk(S)=1$ on the whole of $M$, and the rank 2 bundle $L$ is equal to $S\oplus \wedge^{dim(M)-1} S^{\circ}$.
Otherwise, at 
any point $x$ we have either $S_x=T_xM$  or $L_x=0+\wedge^{dim(M)-1}T^*M$.
In the first case by Cor. \ref{norom} we known that, near $x$, $L$ is the graph of a top form. In the second case  $L$
projects isomorphically onto the second component $\wedge^{dim(M)-1}T^*M$ near $x$, so by Prop. \ref{NPiso} it must be the graph of a $dim(M)$-vector field.
\end{ex}


Finally, we characterize when a regular Lagrangian subbundle is a higher Dirac structure.
\begin{thm}\label{integ}
Let $M$ be a manifold, fix an  integer $ p\ge 1$ and a Lagrangian subbundle $L\subset TM\oplus \wedge^p T^*M$. Assume that  $S:=pr_{TM} (L)$ has constant rank along $M$. Choose a form $\omega\in \Omega^{p+1}(M)$ such that $S$ and $\omega$ describe $L$ as in Cor. \ref{norom}.
 
 Then $L$ is involutive if{f} $S$ is an involutive distribution and $d\omega|_{\wedge^3 S \otimes \wedge^{p-1}TM}=0.$ 
\end{thm}
\begin{proof}
First notice that a $p+1$-form $\omega$ as above always exists, as it can be constructed as in Lemma \ref{ext} choosing a (smooth) distribution $C$ on $M$ complementary to $S$. We use the description of $L$ given in Cor. \ref{norom}.

Assume that $L$ is involutive. By Prop. \ref{algoid}, $S$ is an involutive distribution.
Let $X,Y$ be sections of   $S$. Using $\cL_X\omega=d(\iota_X \omega)+\iota_X d\omega$ we have
$$\oo X+\iota_X \omega, Y+\iota_Y\omega \cc=
[X,Y]+\cL_X(\iota_Y\omega)-\iota_Y (\cL_X\omega)+\iota_Y \iota_X d\omega=
[X,Y]+\iota_{[X,Y]}\omega+\iota_Y \iota_X d \omega.$$
Since this lies in $L$ 
we have $\iota_Y \iota_X d \omega \in \wedge^{p} S^{\circ}$ for all sections $X,Y$ of $S$, which is equivalent to $d\omega|_{\wedge^3 S \otimes \wedge^{p-1}TM}=0.$

Conversely, assume the above two conditions on $S$ and $d\omega$. The above computation shows that for all sections $X,Y$ of $S$, the bracket  $\oo X+\iota_X \omega, Y+\iota_Y\omega \cc$ lies in $L$.
The brackets of $X+\iota_X \omega$ with sections of $\wedge^p S^{\circ}$ lie in $L$ since, by the involutivity of $S$,  locally $\wedge^p S^{\circ}$ admits a frame consisting of $p$-forms $\alpha_i$ which are closed and which hence satisfy $\oo \alpha_i,\cdot \cc=0$. Therefore $L$ is involutive.
\end{proof}

Notice that for $p=1$ (so $d\omega$ is a 3-form) we obtain the familiar statement that a   regular almost Dirac structure $L$ is involutive if{f} $pr_{TM}(L)$ is an involutive distribution whose leaves are endowed with closed 2-forms (see \cite[Thm. 2.3.6]{Cou}).

\section{Equivalence of higher Dirac and Multi-Dirac structures}\label{equivmd}

Recently   Vankerschaver, Yoshimura and Marsden \cite{MultiDirac}	introduced the notion of Multi-Dirac structure. In this section we show that, at least in the regular case, it is equivalent to our notion of higher Dirac structure.
This section does not affect any of the following ones and might be skipped on a first reading.

We recall some definitions from \cite[\S 4]{MultiDirac}. All along we fix an integer $p\ge 1$ and a manifold $M$. In the following the indices $r,s$ range from $1$ to $p$.
Define $$P_r:=\wedge^r TM\oplus \wedge^{p+1-r} T^*M.$$
Define a pairing $P_r\times P_s \to \wedge^{p+1-r-s} T^*M$
by
\begin{equation*} 
\langle\!\langle
 ({Y},  \eta), (\bar{{Y}}, \bar{\eta}) \rangle\!\rangle :=\frac{1}{2}\left( \iota_{\bar{{Y}}}\eta-(-1)^{rs}\iota_{{Y}}\bar{\eta} \right).
\end{equation*}
If $V_s\subset P_s$, then $(V_{s})^{\perp,r}\subset P_r$ is defined by
\begin{equation}\label{perpe} 
(V_{s})^{\perp,r}:= \{  ({Y},  \eta) \in P_{r} :
\langle\!\langle({Y},  \eta)\;,\; V_s \rangle\!\rangle=0  \}.
\end{equation}

\begin{defi}\label{MultiDirac} 
An \emph{almost multi-Dirac structure of degree $p$} on $M$ consists of subbundles
$(D_1, \ldots, D_{p})$, where 
$
D_{r} \subset  P_{r}$ for all $r$,
 satisfying  
\begin{equation} \label{isotropy}
D_{r}=(D_{s})^{\perp, r}
\end{equation} 
for all $r,s$ with $r+s \le p+1$. 
\end{defi}

 
\begin{prop}\label{eqla} Fix a manifold $M$ and an integer $p\ge 1$.
There is a bijection
\begin{align*}
\{\text{almost Multi-Dirac structures of degree }p \} &\cong \{\text{almost Dirac structures }L\text{ of order }p \text{ s.t. } \\
&\;\;\;\;\;\;\;L^{\perp, r} \text{ is a   subbundle for }r=2,\dots,p\}\\
(D_1,\dots,D_p)&\mapsto D_1.
\end{align*}
\end{prop}
The proof of Prop. \ref{eqla} uses the following extension of Cor. \ref{norom}:
\begin{lemma}\label{mnormal} Fix a vector space $T$ and an integer $p\ge 1$.
Let $L$ be a Lagrangian subspace of $T\oplus \wedge^pT^*$, and define
$D_r:=(L)^{\perp, r}$  for $r=1,\dots,p$. Choose   $\omega\in \wedge^{p+1}T^*$ so that $\omega$ and $S:=pr_{T}(L)$ describe $L$ as in Cor. \ref{norom}. Then for all $r$ we have
\begin{equation*}
D_r= \{Y+\iota_Y \omega +\xi: Y\in S\wedge (\wedge^{r-1}T), \xi \in \wedge^{p+1-r} S^{\circ}\}.
\end{equation*}
\end{lemma} 
\begin{proof}
`` $\subset$:'' We first claim that 
 $$pr_{\wedge^r T} (D_r)\subset  S\wedge (\wedge^{r-1}T).$$ If $S=T$ this obvious. 
If $S\neq T$, by Prop. \ref{linalg} we have that $\wedge^p S^{\circ} \subset L$ is non-zero.
As $(Y,\eta)\in D_r$ implies $\iota_Y (\wedge^p S^{\circ})= 0$, we conclude that  $Y\in S\wedge (\wedge^{r-1}T)$.

Let $(Y,\eta)\in D_r$. For all $(X,\alpha)\in L$ we have $\alpha-\iota_{X}\omega \in \wedge^p S^{\circ}$ by Cor. \ref{norom}, and since $Y\in S\wedge (\wedge^{r-1}T)$ we obtain $\iota_{Y}\alpha=\iota_{Y} (\iota_{X}\omega)$. Hence zero equals
\begin{equation}\label{eiota}
\langle\!\langle
 ({Y},  \eta),(X,\alpha)  \rangle\!\rangle = \iota_X \eta -(-1)^r \iota_{Y}\alpha=
 \iota_X \eta -(-1)^r\iota_{Y} (\iota_{X}\omega)=\iota_{X}(\eta-\iota_{Y}\omega),
\end{equation}
that is, $\eta-\iota_{Y}\omega\in \wedge^{p+1-r} S^{\circ}$.
Notice that in the last equality of eq. \eqref{eiota} we used the total skew-symmetry of $\omega$.

`` $\supset$'' follows from eq. \eqref{eiota}.
\end{proof}

\begin{proof}[Proof of Prop. \ref{eqla}] 
The map in the statement  of Prop. \ref{eqla} is  well-defined by eq. \eqref{isotropy} with $r=s=1$.
It is   injective as $D_r=(D_{1})^{\perp, r}$ is determined by $D_1$ for $r=2,\dots,p$, again by eq. \eqref{isotropy}.

 We now show that it is surjective.
Let $L$ be a Lagrangian subbundle of $E^p$, and assume that 
$D_r:=(L)^{\perp, r}$ is a smooth subbundle for $r=1,\dots,p$. We have to show that eq. \eqref{isotropy} holds for all $r,s$ with $r+s\le p+1$.
If $(Y,\eta)\in D_r$ and $(\bar{Y},\bar{\eta})\in D_s$, then $\iota_{Y}\bar{\eta}=\iota_{Y} (\iota_{\bar{Y}}\omega)$ by Lemma \ref{mnormal}, showing $\langle\!\langle
D_r,D_s\rangle\!\rangle=0$ and the 
  inclusion ``$\subset$''.
    
 For the opposite inclusion take  $(Y,\eta)\in (D_{s})^{\perp, r}$ at some point $x\in M$. In particular
 $(Y,\eta)$ is orthogonal to $\wedge^{p+1-s} S^{\circ}_x$ (where $S_x:=pr_{T_xM}L$). The latter does not vanish  by Prop. \ref{linalg} if $S_x\neq T_xM$, and since $r\le p+1-s$ we conclude that $Y\in S_x\wedge (\wedge^{r-1}T_xM)$. If $S_x=T_xM$ the same conclusion holds. A computation analog to eq. \eqref{eiota} implies that for all $(\bar{Y},\bar{\eta})\in D_s$ we have $0=\iota_{\bar{Y}}(\eta-\iota_{Y}\omega)$. As such $\bar{Y}$ span 
$S_x\wedge (\wedge^{s-1}T_xM)$ 
 by    Lemma \ref{mnormal} applied to $D_s$, from $s\le p+1-r$ it follows that $\eta-\iota_{Y}\omega\in \wedge^{p+1-r} S_x^{\circ}$. Hence by  Lemma \ref{mnormal} $(Y,\eta)\in D_r$ .
\end{proof}

In other to introduce the notion of integrability for almost multi-Dirac structures, as in \cite{MultiDirac} define 
$\left[\!\left[\cdot,\cdot  \right]\!\right]_{r,s} \colon \Gamma(P_r)\times \Gamma(P_s) \to \Gamma(P_{r+s-1})$ by 
 \begin{equation*}
\left[\!\left[\left({Y},\eta\right), \left(\bar{{Y}}, \bar{\eta}\right) \right]\!\right]_{r,s}\\
:=
\left( [{Y},\bar{{Y}}], \; \cL_{{Y}}\bar{\eta}-(-1)^{(r-1)(s-1)}\cL_{\bar{{Y}}}\eta+\frac{(-1)}{2}^{r}d
\left( \iota_{\bar{{Y}}}\eta+(-1)^{rs}\iota_{{Y}}\bar{\eta} \right)
\right).
\end{equation*}

\begin{defi} 
An almost Multi-Dirac structure $(D_1, \dots, D_{p})$  is  \emph{integrable}  if 
 \begin{equation} 
  \label{involm}\left[\!\left[D_r, D_s\right]\!\right]_{r,s} \subset D_{r+s-1} 
 \end{equation}
 for all  $r, s$ with  $r+s \le p$. In that case it is a  \emph{Multi-Dirac structure}.
 \end{defi}


We call an almost Multi-Dirac structure $(D_1,\dots,D_p)$  \emph{regular} if $pr_{TM} (D_1)$ has constant rank. By Lemma \ref{mnormal}, this is equivalent to  $pr_{\wedge^r TM} (D_r)$ having constant rank for $r=1,\dots,p$.
Under this regularity assumption, we obtain an equivalence for integrable structures.

\begin{thm}\label{eqint} Fix a manifold $M$ and an integer $p\ge 1$.
The  bijection of Prop. \ref{eqla} restricts to a bijection
\begin{align*}
\{\text{regular Multi-Dirac structures of degree }p \text\} &\cong \{\text{regular Dirac structures of order }p \}
\end{align*}
\end{thm}
\begin{proof} 
If $(D_1,\dots,D_p)$ is a Multi-Dirac structure, by the remark at the end of \cite[\S 4]{MultiDirac}, $D_1$ is involutive w.r.t. the Courant bracket. Therefore it is involutive w.r.t. Dorfman bracket, that is, it is a Dirac structure of order $p$. 

For the converse, notice that if $L$ is  a regular Dirac structure $L$ then 
$L^{\perp, r}$ is always a smooth subbundle by Cor. \ref{mnormal}. 
So let $(D_1,\dots,D_p)$ be a regular almost  Multi-Dirac structure
with the property that $L:=D_1$ is involutive. Choose   $\omega\in \Omega^{p+1}(M)$ so that $(\omega, S:=pr_{TM}(L))$ describe $L$ as in Cor. \ref{norom}. Such a differential form exists by the regularity assumption.
To show that condition \eqref{involm} holds, let  $Y\in \Gamma(S\wedge (\wedge^{r-1}T))$ and $\bar{Y}\in \Gamma(S\wedge (\wedge^{s-1}T))$. 
We have $$\left[\!\left[Y+\iota_Y \omega, \bar{Y}+\iota_{\bar{Y}} \omega\right]\!\right]_{r,s} =\left([Y,  \bar{Y}], \iota_{[Y,  \bar{Y}] }\omega+(-1)^r \iota_Y\iota_{\bar{Y}} d\omega \right),$$ see for instance \cite[Proof of Thm. 4.5]{MultiDirac}.
Now $\iota_Y\iota_{\bar{Y}} d\omega \in \Gamma(\wedge^{p+2-r-s}S^{\circ})$ by Thm.  \ref{integ}, so the above lies in $D_{r+s-1}$ by Lemma \ref{mnormal}.
Further, the involutivity of $S$   implies that  locally $\wedge^{p+1-s} S^{\circ}$ admits a frame consisting of closed forms $\alpha_i$. For any choice of functions $f_i$ we have
 $$\left[\!\left[Y+\iota_Y \omega, f_i \alpha_i \right]\!\right]_{r,s}
=\cL_Y (f_i \alpha_i) + (-1)^{r(s+1)} d \iota_Y(f_i \alpha_i)
=\iota_{Y}(df_i\wedge \alpha_i),$$ which lies in  $\Gamma(\wedge^{p+2-r-s}S^{\circ})$
since $Y\in \Gamma(S\wedge (\wedge^{r-1}T))$ and $\alpha_i \in \Gamma(\wedge^{p+1-s} S^{\circ})$.
\end{proof}



Finally, we comment on how our definition of higher Dirac structure differs from 
Hagiwara's Nambu-Dirac structures \cite{Hagi}, which   also are an extension of  
Courant's notion of Dirac structure. 
\begin{remark}\label{Hagi}
A  \emph{Nambu-Dirac structure} on a manifold $M$ \cite[Def. 3.1, Def. 3.7]{Hagi} is an involutive subbundle $L\subset E^p$ satisfying  
\begin{align}\label{Hagiiso}
&\la X_1+\alpha_1,X_2+\alpha_2\ra|_{\wedge^{p-1}(pr_{TM}(L))}=0,\\
\label{hismax}
&\wedge^{p}(pr_{TM}(L))=pr_{\wedge^{p}TM}L^{\perp,p},
\end{align}
 where $L^{\perp,p}\subset \wedge^{p}TM\oplus T^*M$ is defined as in eq. \eqref{perpe}.   When $p=1$, Nambu-Dirac structures agree with Dirac structures. Graphs of closed forms and of Nambu-Poisson multivector fields are Nambu-Dirac structures. 

Our isotropicity condition \eqref{isot} is clearly stronger than \eqref{Hagiiso}. Nevertheless, higher Dirac structures are usually \emph{not} Nambu-Dirac structures, for the former satisfy
$$pr_{TM}(L)\wedge(\wedge^{p-1}TM)= pr_{\wedge^{p}TM}L^{\perp,p}$$ by Lemma \ref{mnormal}, and hence usually do not satisfy \eqref{hismax}. 
 A concrete instance is given by the 3-dimensional Lagrangian  subspace $L\subset T\oplus \wedge^2 T^*$ given as in Cor. \ref{norom} by $T=\RR^4$, $S$ equal to the plane
$\{x_3=x_4=0\}$ and $\omega=dx_1\wedge dx_2\wedge dx_3$.\end{remark}

\section{Review:  $L_{\infty}$-algebras}\label{Linfty}
 
In this section we review briefly the notion of $L_{\infty}$-algebra, which generalizes  Lie algebras    
and was introduced by Stasheff \cite{LadaStasheff} in the 1990s. 
We will follow the conventions of Lada-Markl\footnote{Except that on graded vector spaces we take the grading inverse to theirs.}
 \cite[\S2,\S5]{LadaMarkl}.

Recall that a \emph{graded vector space} is just a (finite dimensional, real) vector space $V=\oplus_{i\in \ZZ}V_i$ with a direct sum decomposition into subspaces. An element  
of $V_i$ is said to have  degree $i$, and we 
denote its degree by $|\cdot|$.

For any $n\ge 1$, $V^{\otimes n}$ is a graded vector space, and the symmetric group 
acts on it 
by the 
so-called odd representation: the transposition 
of the $k$-th and $(k+1)$-th element  acts by
$$v_1\otimes\dots\otimes v_n \mapsto -(-1)^{|v_k||v_{k+1}|} v_1\otimes\dots\otimes 
v_{k+1}\otimes v_k \otimes\dots\otimes v_n.$$
The \emph{$n$-th graded exterior product of $V$} is the graded vector space
 $\wedge^n V$,  consisting of elements  of $V^{\otimes n}$ which are fixed by the odd representation of the symmetric group.
\begin{defi}\label{lidef}
An  \emph{$L_{\infty}$-algebra}  is a  graded vector space $V=\bigoplus_{i\in \ZZ}V_i$ endowed with a sequence of multi-brackets ($n\ge 1$)
\begin{equation*}
l_n \colon \wedge ^n V \to V
\end{equation*}
of degree $2-n$, satisfying the following quadratic relations  for each $n\ge 1$:
  \begin{align}\label{lijac} 
\sum_{i+j=n+1}\sum_{\sigma\in Sh(i,n-i)}\chi(\sigma)(-1)^{i(j-1)}l_j(l_i(v_{\sigma(1)},\dots,v_{\sigma(i)}),
v_{\sigma(i+1)},\dots,v_{\sigma(n)})=0.
\end{align}
Here $Sh(i,n-i)$ denotes the set of $(i,n-i)$-unshuffles, that is, permutations 
preserving the order of the first $i$ elements and the order of the last $n-i$ elements.
The sign $\chi(\sigma)$ is given by the action of $\sigma$ on 
$v_1\otimes\dots\otimes v_n$ in the  odd representation.
\end{defi} 

\begin{remark}
1) The quadratic relations imply that the unary bracket $l_1$ squares to zero, so  $(V,l_1)$ is a chain complex of vector spaces. Hence  $L_{\infty}$-algebras can be viewed as chain complexes with the extra data given by the multi-brackets  
$l_n$ for $n\ge 2$.

2) When $V$ is concentrated in degree $0$,
(i.e., only $V_0$ is non-trivial) then  $\wedge^n V$ is the  usual $n$-th exterior product of $V$, and is concentrated in degree zero. Hence
 by degree reasons only the  binary bracket $[\cdot,\cdot]_2$ is non-zero, and the   quadratic relations  are simply the Jacobi identity, so we recover the notion of Lie algebra.
\end{remark}

For any $p\ge 1$, we use the term \emph{Lie $p$-algebra} to denote an 
$L_{\infty}$-algebra whose underlying graded vector space is concentrated in   degrees $-p+1,\cdots,0$. Notice that
by degree reasons 
only the multi-brackets $l_1,\cdots,l_{p+1}$ can be non-zero.
In particular, a Lie 2-algebra consists of a graded vector space $V$
concentrated in degrees $-1$ and $0$, together with
maps
\begin{align*}
d:=&l_1 \colon V\to V\\        
[\cdot,\cdot]:=&l_2 \colon \wedge ^2 V \to V\\
J:=&l_3 \colon \wedge ^3 V \to V
\end{align*}
of degrees $1$,$0$ and $-1$ respectively, subject to the  quadratic relations.\\

An \emph{$L_{\infty}$-morphism} $\phi \colon V \rightsquigarrow V'$  between $L_{\infty}$-algebras is a sequence of maps ($n\ge 1$)
\begin{equation*}
\phi_n \colon \wedge ^n V \to V'
\end{equation*}
of degree $1-n$, satisfying certain relations, which can be found in
 \cite[Def. 5.2]{LadaMarkl} in the case when $V'$ has only the unary and binary bracket.
The first of these relations says that $\phi_1 \colon V \to V'$ must preserve the differentials  (unary brackets). We spell out the definition when 
$V$ and $V'$ are Lie 2-algebras.

 
\begin{defi}\label{defmor} Let $(V,d,[\cdot,\cdot],J)$ and $(V',d',[\cdot,\cdot]',J')
$ be Lie 2-algebras.
A \emph{morphism}   $\phi \colon V \rightsquigarrow V'$  consists  of  
linear maps 
\begin{align*}
\phi_0  &\colon  V_0 \to V_0\\
\phi_{1}  &\colon  V_{-1} \to V_{-1}\\
\phi_2 & \colon\wedge^2 V_0 \to V_{-1}        
\end{align*}
such that 
 \begin{align} \label{chainmap}d' \circ\phi_{1}&= \phi_{0}\circ d,\\
\label{failjac}d' (\phi_2(x,y))&=\phi_0[x,y]-[\phi_0(x),\phi_0(y)]'
\;\;\;\text{ for all } x,y \in V_0,\\
\label{failjacnew}\phi_2(df,y)&=\phi_1[f,y]-[\phi_1(f),\phi_0(y)]'
\;\;\;\text{ for all } f\in V_{-1},y \in V_0,
\end{align}
and  for all $x,y,z \in V_0$:
 \begin{align}\label{eight}
&\phi_0(J(x,y,z))-J'(\phi_0(x),\phi_0(y),\phi_0(z))=\\
&\phi_2(x,[y,z])
-\phi_2(y,[x,z])+\phi_2(z,[x,y]) \nonumber \\
+&[\phi_0(x),\phi_2(y,z)]'-[\phi_0(y),\phi_2(x,z)]'
+[\phi_0(z),\phi_2(x,y)]'  \nonumber.
\end{align}
\end{defi}

\section{$L_{\infty}$-algebras from higher analogues of Dirac structures}\label{obs}

Courant \cite[\S 2.5]{Cou} associated to every Dirac structure on $M$ a subset of $C^{\infty}(M)$, which we refer to as Hamiltonian functions or observables.
Usually the Hamiltonian vector field associated to such a function is not unique. Nevertheless, the set of  Hamiltonian functions is endowed with
a Poisson algebra structure (a Lie bracket compatible with the product of functions).
Baez, Rogers and Hoffnung associate to a $p$-plectic form a set of Hamiltonian $p-1$-forms and endow it with a bracket \cite[\S 3]{BHR}. Rogers shows that  the bracket can be extended to obtain a Lie $p$-algebra  \cite[Thm. 5.2]{RogersL}.
In this section we mimic Courant's definition of the bracket and extend Roger's  results to arbitrary isotropic involutive subbundles.

Let  $p\ge 1$ and let $L$ be an isotropic, involutive subbundle of $E^p=TM\oplus \wedge^pT^*M$.

\begin{defi}\label{ham} 
A $(p-1)$-form $\alpha \in \Omega^{p-1}(M)$ is called \emph{Hamiltonian} if there exists a smooth vector field $X_{\alpha}$ such that $X_{\alpha}+d\alpha \in \Gamma(L)$. We denote the set of Hamiltonian forms by $\Omega^{p-1}_{ham}(M,L)$.
We refer to $X_{\alpha}$ as \emph{a Hamiltonian vector field} of $\alpha$.
\end{defi}

\begin{remark}\label{ann}
a) Hamiltonian vector fields are unique only up to smooth sections of $L\cap (TM\oplus 0)$.

b) For all $X\in L_x\cap (T_xM\oplus 0)$ and for all  $\eta \in pr_{\wedge^pT^*M} L_x$ , 
$$\iota_X \eta=0.$$ Here $x\in M$ and $pr_{\wedge^pT^*M}$ denotes the projection of $E^p_x$ onto the second component. The above property follows from the fact that there exists $Y\in T_xM$ with $Y+\eta \in L_x$, so
$\iota_X \eta = \la X+0,Y+\eta \ra=0$ by the isotropicity of $L$.
\end{remark}

\begin{defi}\label{brs} We define a bracket
$\{\cdot,\cdot\}$ on $\Omega^{p-1}_{ham}(M,L)$ by
\begin{equation*}
\{\alpha,\beta\}:= \iota_{X_{\alpha}} d\beta, 
\end{equation*}
where $X_{\alpha}$ is any Hamiltonian vector field for  $\alpha$.
\end{defi}

\begin{lemma}\label{bracket}
The bracket $\{\cdot,\cdot\}$ is well-defined and skew-symmetric.
It does not satisfy the Jacobi identity, but rather
\begin{equation*}
\{\alpha,\{\beta,\gamma\}\} +c.p.=-d (\iota_{X_{\alpha}}\{\beta,\gamma\})
\end{equation*}
where ``$c.p.$'' denotes cyclic permutations. 
\end{lemma}
 
\begin{proof}
The bracket is well-defined: by Remark \ref{ann}
the ambiguity in
  the choice of $X_{\alpha}$ is
  a section $X$ of $L\cap (TM\oplus 0)$ and  $\iota_X d \beta =0$.
Using $\cL_Y=\iota_Yd+d\iota_Y$ one computes  
\begin{equation}\label{closed}
\oo X_{\alpha}+d\alpha, X_{\beta}+d\beta\cc= [X_{\alpha},X_{\beta}]+ d \{ \alpha, \beta\}.
\end{equation}
 Hence $[X_{\alpha}, X_{\beta}]$ is a Hamiltonian vector field for $\{\alpha,\beta\}$, showing that $\Omega^{p-1}_{ham}(M,L)$ is closed under $\{\cdot,\cdot\}$.
The bracket is skew symmetric because   $$0=\la X_{\alpha}+d\alpha, X_{\beta} +d\beta \ra=
\{\alpha,\beta\}+\{\beta,\alpha\}.$$

To compute the Jacobiator of $\{\cdot,\cdot\}$ we proceed as in\footnote{There the case $p=1$ is treated, and the term $\iota_{X_{\alpha}}\{\beta,\gamma\} $ vanishes by degree reasons.} \cite[Prop. 2.5.3]{Cou}. Since $L$ is isotropic and involutive we have
\begin{align} 
0&=\la \oo X_{\alpha}+d\alpha\;,\;X_{\beta}+d\beta \cc,X_{\gamma}+d\gamma \ra  \nonumber \\
&=\la [X_{\alpha},X_{\beta}]+d\{\alpha,\beta\}\;,\;X_{\gamma}+d\gamma \ra  \nonumber \\
  &= \iota_{[X_{\alpha},X_{\beta}]}d\gamma+ \iota_{X_{\gamma}}d\{\alpha,\beta\}  \nonumber \\
&= \left(\{ \alpha,\{\beta,\gamma\}\} +c.p.\right)+d \;(\iota_{X_{\alpha}}\{\beta,\gamma\}). \nonumber
\end{align}
Here the second equality   uses eq. \eqref{closed} and the last equality uses $\iota_{[Y,Z]}=[\cL_Y,\iota_Z]$.
\end{proof}

\begin{remark}\label{graph-o}
Given a $p$-plectic form $\omega$,
Cantrijn, Ibort and de Le\'on \cite[\S 4]{Hammulti} 
define the   space of 
 Hamiltonian $(p-1)$-forms  $\alpha$  by the requirement that $d\alpha=-\iota_{X_{\alpha}}\omega$
for a (necessarily unique) vector field $X_{\alpha}$ on $M$, and define 
the \emph{semi-bracket}  $\{\alpha,\beta\}_s$ by $\iota_{X_{\beta}}\iota_{X_{\alpha}}\omega$. 
These notions
coincide with our Def. \ref{ham} and Def. \ref{brs} applied to $graph(\omega):=\{X-\iota_X \omega: X\in TM\}\subset E^{p}$. 
\end{remark}

\begin{remark} Given an $p$-plectic form,
in   \cite[Def. 3.3]{BHR} the \emph{hemi-bracket} of $\alpha,\beta \in
\Omega^{p-1}_{ham}(M,graph(\omega))$ is also defined, by the  formula $\cL_{X_{\alpha}}\beta$. 
This notion does not extend to the setting of arbitrary isotropic subbundles of $E^p$, since in that setting 
the Hamiltonian vector field $X_{\alpha}$ is not longer unique and the above expression depends on it.

For instance, take $M=\RR^4$, consider the closed $3$-form $\theta=dx_1\wedge dx_2 \wedge dx_3$. By Prop. \ref{closedform}, $L=\{X-\iota_X \theta: X\in TM\}$ 
 is a isotropic, involutive subbundle of $E^2$. Both $\pd{x_4}\in \Gamma(L\cap TM)$ and the zero vector field are Hamiltonian vector fields for $\alpha=0$, and 
 the hemi-bracket of $\alpha$ with $\beta=x_1dx_4+x_4 dx_1$ is not well-defined since
$$\cL_{\pd{x_4}}\beta=dx_1\neq 0=\cL_{0}\beta.$$ \end{remark}

 Rogers  \cite[Thm. 5.2]{RogersL} shows that for every $p$-plectic manifold there is an associated $L_{\infty}$-algebra of observables. The statement and the proof generalize in a straightforward way to arbitrary isotropic, involutive subbundle of $E^p=TM\oplus \wedge^p T^*M$.
   
\begin{thm}\label{Liep}
Let $p\ge 1$ and $L$ be a isotropic, involutive subbundle of $E^p=TM\oplus \wedge^p T^*M$. Then the  complex concentrated in degrees $-p+1,\dots,0$
$$
C^{\infty}(M) \overset{d}{\rightarrow}\dots \overset{d}{\rightarrow}\Omega^{p-2}(M)\overset{d}{\rightarrow}\Omega^{p-1}_{ham}(M,L)$$ has a Lie $p$-algebra structure. 
The only non-vanishing multibrackets are given by the de Rham differential  on $\Omega^{\le p-2}(M)$ and, for $k=2,\dots,p+1$, by 
$$
l_k(\alpha_1,\dots,\alpha_k)=
\;\;\;\epsilon(k)\iota_{X_{\alpha_k}}\dots \iota_{X_{\alpha_3}}\{\alpha_1,\alpha_2\}
$$
where $\alpha_1,\dots,\alpha_k \in \Omega^{p-1}_{ham}(M,L)$ and  $\epsilon(k)=(-1)^{\frac{k}{2}+1}$ if $k$ is even, 
$\epsilon(k)=(-1)^{\frac{k-1}{2}}$ if $k$ is odd.
\end{thm}

%

\begin{proof}
The expressions for the multibrackets are totally skew-symmetric, as a consequence of the fact that  $\{\cdot,\cdot\}$ is skew-symmetric. This and the fact that $\{\cdot,\cdot\}$ is independent of the choice of Hamiltonian vector fields imply that the multibrackets are well-defined.
Clearly   $l_k$ has degree $2-k$. 

Now we check the $L_{\infty}$ relations   \eqref{lijac}. For $n=1$ the relation holds due to $d^2=0$. Now consider the relation \eqref{lijac} for a fixed $n\ge 2$, and let $\alpha_1,\dots,\alpha_n$ be homogeneous elements of the above complex. We will use repeatedly the fact that, for $k\ge2$,  the $k$-multibracket  vanishes when one of its entries is of negative degree.
For $j\in \{2,\dots,n-2\}$ (so $i\ge 3$), we have
$$l_j(l_i(\alpha_{1}, \dots,\alpha_{i}),\alpha_{i+1}, \dots, \alpha_{n} )=0,$$ as a consequence of the fact that  $k$-multibrackets for $k\ge3$ take values in negative
 degrees. 
For $j=n$ we have  
$$l_n(l_1(\alpha_{1}),  \alpha_{2}, \dots, \alpha_{n} )=0:$$
if $|\alpha_{1}|=0$ then $l_1(\alpha_{1})$ vanishes, otherwise 
$l_1(\alpha_{1})=d\alpha_{1}$ and its Hamiltonian vector field vanishes.
 
 We are left with the summands of  \eqref{lijac} with   $j=1$ and $j=n-1$.
When $n=2$ we have just one summand $l_1(l_2(  \alpha_{\sigma(1)},\alpha_{\sigma(2)}))$ 
which vanishes by degree reasons.
 For $n\ge 3$ it is enough to assume that  all the $\alpha_i$'s have degree zero. We have 
\begin{align*}
d(l_n(\alpha_{ 1}, \dots,\alpha_{ n}))+
\sum_{\sigma
\in Sh(2,n-2)} \chi(\sigma)   l_{n-1}(\{\alpha_{\sigma(1)},\alpha_{\sigma(2)}\},\alpha_{\sigma(3)} \dots, \alpha_{\sigma(n)} ).        
\end{align*}
Writing out explicitly the unshuffles in $Sh(2,n-2)$ and the multibrackets we obtain
\begin{align*}
\epsilon(n)&d(\iota_{X_{\alpha_{n}}}\dots\iota_{X_{\alpha_{3}}}\{\alpha_{1},\alpha_{2}\})\\
+ \epsilon(n-1)&\Big[
  \sum_{2\le i<j\le n }(-1)^{i+j-1}\iota_{X_{\alpha_{n}}}\dots  \widehat{\iota_{X_{\alpha_{j}}}}\dots
 \widehat{\iota_{X_{\alpha_{i}}}}
\dots
\iota_{X_{\alpha_{2}}}\{\{\alpha_{i},\alpha_{j}\},\alpha_{1}\} \\\
&+ \sum_{3\le j\le n }(-1)^{j}\iota_{X_{\alpha_{n}}}\dots  \widehat{\iota_{X_{\alpha_{j}}}} 
\dots
\iota_{X_{\alpha_{3}}}\{\{\alpha_{1},\alpha_{j}\},\alpha_{2}\} \\\
&+ \iota_{X_{\alpha_n}}\dots   
\dots
\iota_{X_{\alpha_{4}}}\{\{\alpha_{1},\alpha_{2}\},\alpha_{3}\} \
\Big].
\end{align*}
By Lemma \ref{pain} 
 we conclude that the above expression vanishes.
 \end{proof}

The following Lemma, needed in the proof of Thm. \ref{Liep},
 extends \cite[Lemma 3.7]{RogersL}.
   
\begin{lemma}\label{pain}
Let $p\ge 1$ and $L$ be a isotropic, involutive subbundle of $E^p=TM\oplus \wedge^p T^*M$. Then for any $n\ge 3$, and for all $\alpha_1,\dots,\alpha_n\in \Omega^{p-1}_{ham}(M,L)$ we have
\begin{align*}
d(\iota_{X_{\alpha_n}}\dots\iota_{X_{\alpha_3}}\{\alpha_{1},\alpha_{2}\})=
 (-1)^{n+1}&\Big[
  \sum_{2\le i<j\le n }(-1)^{i+j-1}\iota_{X_{\alpha_n}}\dots  \widehat{\iota_{X_{\alpha_j}}}\dots
 \widehat{\iota_{X_{\alpha_i}}}
\dots
\iota_{X_{\alpha_2}}\{\{\alpha_{i},\alpha_{j}\},\alpha_{1}\} \\\
&+ \sum_{3\le j\le n }(-1)^{j}\iota_{X_{\alpha_n}}\dots  \widehat{\iota_{X_{\alpha_j}}} 
\dots
\iota_{X_{\alpha_3}}\{\{\alpha_{1},\alpha_{j}\},\alpha_{2}\} \\\
&+ \iota_{X_{\alpha_n}}\dots   
\dots
\iota_{X_{\alpha_4}}\{\{\alpha_{1},\alpha_{2}\},\alpha_{3}\} \
\Big].
\end{align*}
\end{lemma}
\begin{proof}
We proceed by induction on $n$. For $n=3$ the statement holds by Lemma \ref{bracket}. So let $n>3$. To shorten the notation, denote $A:=\iota_{X_{\alpha_{n-1}}}\dots\iota_{X_{\alpha_3}}\{\alpha_{1},\alpha_{2}\}$. Then we have
\begin{align}\label{lhslie}
d(\iota_{X_{\alpha_n}}\dots\iota_{X_{\alpha_3}}\{\alpha_{1},\alpha_{2}\})=d(\iota_{X_{\alpha_n}} A)=\cL_{X_{\alpha_n}}A-\iota_{X_{\alpha_n}}dA.
\end{align}
The first term on the r.h.s. of \eqref{lhslie}
becomes 
\begin{align*} 
&\cL_{X_{\alpha_n}}(\iota_{X_{\alpha_{3}}\wedge \dots\wedge {X_{\alpha_{n-1}}}}\{\alpha_{1},\alpha_{2}\})\\
=& \sum_{i=3}^{n-1}(-1)^{i+1}\iota_{X_{\alpha_{n-1}}}\dots \widehat{\iota_{X_{\alpha_i}}}\dots\iota_{X_{\alpha_3}}
\iota_{[X_{\alpha_n},X_{\alpha_i}]
}
\{\alpha_{1},\alpha_{2}\}
+\iota_{X_{\alpha_{n-1}}}\dots \iota_{X_{\alpha_3}}
\cL_{X_{\alpha_n}}\{\alpha_{1},\alpha_{2}\}
\\
=&\sum_{i=3}^{n-1}(-1)^{i+1}\iota_{X_{\alpha_{n-1}}}\dots \widehat{\iota_{X_{\alpha_i}}}\dots 
\iota_{X_{\alpha_2}} \{\{\alpha_n,\alpha_i\},\alpha_{1}\}
+
\iota_{X_{\alpha_{n-1}}}\dots \iota_{X_{\alpha_3}}(
\{\{\alpha_2,\alpha_n\},\alpha_{1}\}-\{\{\alpha_1,\alpha_n\},\alpha_{2}\}
)\\
=&\sum_{i=2}^{n-1}(-1)^{i}\iota_{X_{\alpha_{n-1}}}\dots \widehat{\iota_{X_{\alpha_i}}}\dots\iota_{X_{\alpha_2}} \{\{\alpha_i,\alpha_n\},\alpha_{1}\}
- \iota_{X_{\alpha_{n-1}}}\dots \iota_{X_{\alpha_3}}\{\{\alpha_1,\alpha_n\},\alpha_{2}\}.\end{align*}
Here in the second equality  we used  $[X_{\alpha_n},X_{\alpha_i}]=X_{\{\alpha_n,\alpha_i\}}$ (see the proof of Lemma \ref{bracket})
 and $$\iota_{X_{\{\alpha_n,\alpha_i\}}}\{\alpha_{1},\alpha_{2}\}=-\iota_{X_{\{\alpha_{n},\alpha_{i}\}}}\iota_{X_{\alpha_2}} d \alpha_1=\iota_{X_{\alpha_2}} \{\{\alpha_n,\alpha_i\},\alpha_{1}\},$$
as well as Cartan's formula for the Lie derivative and  Lemma \ref{bracket}.

The second term on the r.h.s. of \eqref{lhslie} can be developed using the induction hypothesis. The resulting expression for the l.h.s. of eq. \eqref{lhslie} 
  is easily seen to agree with the one in the statement of this lemma. 
\end{proof} 

\begin{remark} The observables associated by Thm. \ref{Liep} to the zero $p+1$-form on $M$  are given by the abelian Lie algebra $\RR$ for $p=1$ and to the complex $C^{\infty}(M) \overset{d}{\rightarrow}  \Omega^{1}_{closed}(M)$ (with vanishing higher brackets) for $p=2$. 
It is a curious coincidence that they agree with the central extensions of observables of $p$-plectic structures  given in \cite[Prop. 9.4]{RogersPre} for $p=1$ and $2$ respectively.
\end{remark}


A closed 2-form $B$ on $M$ induce an automorphism of the Courant algebroid $TM\oplus T^*M$ by gauge transformations
(see  \S \ref{intro}), and therefore acts on the set of Dirac structures. For instance, the Dirac structure $TM\oplus \{0\}$ is mapped to the graph of $B$.
The Poisson algebras of observables of  these two Dirac structures are not isomorphic (unless $B=0$). 

Similarly, for $p\ge1$, gauge transformations of $E^p$ by closed $p+1$-forms usually do not
 induce an isomorphism of the Lie $p$-algebra of observables.  
We display a quite trivial operation which, on the other hand, does have this property.

\begin{lemma}\label{lambda} Let $\lambda \in \RR-\{0\}$ and consider
\begin{align*}
m_{\lambda} \colon\;\;\;\;\; E^p &\to E^p\\ X+\eta& \mapsto X+\lambda \eta  
\end{align*}
Let $L\subset E^p$ be an involutive isotropic subbundle.
Then $m_{\lambda}(L)$ is also an involutive isotropic subbundle, and the Lie $p$-algebras of observables of $L$ and $m_{\lambda}(L)$ are isomorphic.
\end{lemma}
\begin{proof}
$m_{\lambda}$ is an automorphism of the Dorfman bracket $\oo \cdot,\cdot \cc$ and 
$\langle m_{\lambda}\cdot,m_{\lambda}\cdot \ra=\lambda \langle  \cdot,\cdot\ra$. Hence 
 $m_{\lambda}(L)$ is also involutive and isotropic. 
  
 We consider the Lie $p$-algebras of observables associated to $L$ and $m_{\lambda}(L)$ respectively, as in Thm. \ref{Liep}. We denote them by $\cO^L$ and  $\cO^{m_{\lambda}(L)}$ respectively.   The underlying complexes coincide, both being
$$ C^{\infty}(M) \overset{d}{\rightarrow}\Omega^1(M)\overset{d}{\rightarrow}\dots\overset{d}{\rightarrow}\Omega^{p-1}_{ham}(M,L).$$

 Notice that if $\alpha\in \Omega^{p-1}_{ham}(M,L)$ has Hamiltonian vector field
 $X_{\alpha}^L$, then $\lambda \alpha$ is a Hamiltonian $(p-1)$-form for  $m_{\lambda}(L)$,  and   $X_{\alpha}^L$ itself is a    Hamiltonian vector field for it. Hence
from   Thm. \ref{Liep} it is clear that the unary map 
given by multiplication by $\lambda$
$$\phi \colon (\beta_0,\dots,\beta_{p-1}) \mapsto (\lambda \beta_0,\dots, \lambda \beta_{p-1})$$  intertwines the multibrackets of $\cO^L$ and $\cO^{m_{\lambda}(L)}$, where $\beta_i\in \Omega^i(M)$ for $i< p-1$ and $\beta_{p-1}\in \Omega^{p-1}_{ham}(M,L)$.
 Therefore, 
setting the higher maps   to zero,
we obtain a strict morphism \cite[\S 7]{MarlkDoubek}
of  Lie $p$-algebras, which clearly is an isomorphism.   
\end{proof} 

As an application of Lemma \ref{lambda} we show that 
to any compact, connected, orientable $p+1$-dimensional manifold  ($p\ge 1$) 
there is an associated  Lie $p$-algebra.
 A dual version of this Lie $p$-algebra appeared in  \cite[Thm. 6.1]{RogerVol}.

\begin{cor}\label{M3}
Let $M$ be a compact, connected, orientable $p+1$-dimensional manifold. For any  volume form $\omega$ consider the Lie $p$-algebra associated 
to $graph(\omega)$ by   Thm. \ref{Liep}, whose underlying complex is
$$ C^{\infty}(M) \overset{d}{\rightarrow}\Omega^1(M)\overset{d}{\rightarrow}\dots\overset{d}{\rightarrow}\Omega^{p-1}(M).$$
 (Notice that all $p-1$-forms are Hamiltonian). Its isomorphism class is independent of the choice of $\omega$, and therefore depends only on the manifold $M$.
\end{cor}
\begin{proof}
Let $\omega_0$ and $\omega_1$ be two volume forms on $M$. They define non-zero cohomology classes in 
$H^{p+1}(M,\RR)=\RR$, so there is a (unique) $\lambda \in \RR-\{0\}$ such that $[\omega_1]=
\lambda [\omega_0]$. By Moser's theorem \cite{Moser}
there is a diffeomorphism $\psi$ of $M$ such that $\psi^*(\omega_1)=\lambda \omega_0$.
This explains the first isomorphism in 
$$\text{Lie $p$-algebra of }\omega_1 \;\cong\; \text{Lie $p$-algebra of }\lambda \omega_0\;\cong \;\text{Lie $p$-algebra of }\omega_0,$$
whereas the second one holds by Lemma \ref{lambda}.
\end{proof}


\section{Relations to $L_{\infty}$-algebras arising from split Courant algebroids}\label{liM}

In this section we construct an $L_{\infty}$-morphism from a Lie algebra associated to $E^0$ with the $\sigma$-twisted bracket, where  $\sigma$ is a closed  2-form, to a Lie 2-algebra associated to $E^1$ with the untwisted Courant bracket (in other words, the Courant bracket twisted by $d\sigma =0$).

We consider again $E^p:=TM\oplus \wedge^pT^*M$.
For $p=0$ we have $E^0=TM \oplus \RR$. Fix a closed 2-form $\sigma \in \Omega^2_{closed}(M)$.
Then $\Gamma(E^0)$ with the 
 $\sigma$-twisted Dorfman bracket 
\begin{equation*}
[X+f,Y+g]_{\sigma}=[X,Y]+(X(g)-Y(f))+\sigma(X,Y)
\end{equation*}
is an honest Lie algebra. (See
\cite[\S 3.8]{Gu}, where a geometric interpretation in terms of circle bundles is given too.)

For $p=1$ we have the (untwisted) Courant algebroid $E^1=TM\oplus T^*M$. Roytenberg and Weinstein \cite{rw} associated to it an  $L_{\infty}$-algebra. In the version given in  \cite[Thm. 4.4]{RogersCou} the underlying complex is
\begin{equation}\label{fctE}
C^{\infty}(M) \overset{d}{\rightarrow}\Gamma(E^1)
\end{equation}
where $d$ is  the de Rham differential. The   binary bracket $[\cdot,\cdot]'$ is given by 
  the \emph{Courant} bracket $\oo \cdot,\cdot \cc_{Cou}$   on $\Gamma(E^1)$ and by
 $$[e,f]'=-[f,e]':=\frac{1}{2}\la e, df \ra$$ for $e\in \Gamma(E^1)$ and $f\in C^{\infty}(M)$.
The trinary bracket $J'$ is given by
 $$J'(e_1,e_2,e_3)=-\frac{1}{6}\left( \la \oo e_1,e_2 \cc_{Cou},e_3\ra + c.p.\right)$$ for elements of $\Gamma(E^1)$, where ``c.p.'' denotes cyclic permutation.
All other brackets vanish. 

We show that there is a canonical morphism between these two  Lie 2-algebras:
\begin{thm}\label{mor01} Let $M$ be a manifold and  $\sigma \in \Omega^2_{closed}(M)$.
There is a canonical  morphism of Lie 2-algebras
\begin{equation}
\phi \colon \left(\Gamma(E^0), [\cdot,\cdot]_{\sigma}\right)   \;\;\rightsquigarrow  \;\; \left(C^{\infty}(M) \overset{d}{\rightarrow}\Gamma(E^1),[\cdot,\cdot]',J'\right)
\end{equation}
 given by
 \begin{align*}
&\phi_0 \colon \Gamma(E^0) \to \Gamma(E^1), \;\;\;\;\;\;\;\;\;\;\;\; (X,f)\mapsto (X,df)\\
 &\phi_2 \colon \wedge^2\Gamma(E^0) \to C^{\infty}(M) , \;\;\; (X,f),(Y,g)\mapsto \frac{1}{2}\big(X(g)-Y(f)\big)+\sigma(X,Y).
\end{align*}
\end{thm} 
\begin{proof}
We check that the conditions of Def. \ref{defmor} are satisfied.
Eq. \eqref{chainmap} is satisfied because $ \Gamma(E^0)$ is concentrated in degree zero.

Eq. \eqref{failjac} is satisfied because for any $X+f,Y+g\in \Gamma(E^0)$ we have
\begin{align*}
& \phi_0\Big[X+f,Y+g\Big]_{\sigma}-\Big[\!\Big[ \phi_0 (X+f),\phi_0 (Y+g)\Big]\!\Big]_{Cou}\\
=&\Big([X,Y]+d\big(X(g)-Y(f)+\sigma(X,Y)\big)\Big)-
\Big	([X,Y]+\frac{1}{2}d(X(g)-Y(f))\Big)
\\
=&d\Big(\phi_2(X+f,Y+g)\Big).
\end{align*}
Eq. \eqref{failjacnew} is satisfied because $\Gamma(E^{\circ})$ is concentrated in degree zero.

We are left with checking eq. \eqref{eight}. 
Let $X+f,Y+g,Z+h\in \Gamma(E^1)$. We want to show that
\begin{align}\label{consist}
-J'(X+df,Y+dg,Z+dh)\;\overset{!}{=}\;&\phi_2\Big(X+f,[Y,Z]
+Y(h)-Z(g)+\sigma(Y,Z)\Big)+c.p.\\+& [ X+df,\phi_2(Y+g,Z+h)]'+c.p. \nonumber
\end{align}
where as usual ``$c.p.$'' denotes cyclic permutation.
The l.h.s. of  eq. \eqref{consist} is equal to
\begin{align*}
&\frac{1}{6}\Big( \la\oo X+df,Y+dg\cc_{Cou},Z+dh\ra \Big)+ c.p. \\
= & \frac{1}{6}\Big( [X,Y](h)+\frac{1}{2}Z(X(g))- \frac{1}{2}Z(Y(f)) \Big)+c.p. \\
=& \frac{1}{4}[X,Y](h) + c.p.
\end{align*}
The r.h.s. is equal to
\begin{align*}
&\frac{1}{2}\Big(X\big(Y(h)-Z(g)+\sigma(Y,Z)\big)-[Y,Z](f)\Big)+\sigma(X,[Y,Z]) + c.p.  \\
&+\frac{1}{2}X\Big(\frac{1}{2}(Y(h)-Z(g))+\sigma(Y,Z)\Big)+c.p. \\
 =&\frac{3}{4}X\Big(Y(h)-Z(g)\Big)-\frac{1}{2}[Y,Z](f)+c.p.\\
 &+\sigma(X,[Y,Z])+X(\sigma(Y,Z))+ c.p. \\
=&  \frac{1}{4}[X,Y](h) + c.p.\\
&+d\sigma(X,Y,Z).
\end{align*}
Since $\sigma$ is a closed form, we conclude that  eq. \eqref{consist} is satisfied.
\end{proof}

\section{$L_{\infty}$-algebras from higher analogues of split Courant algebroids}\label{ez}

In this section we apply Getzler's recent contruction \cite{GetzlerHigherDer} to obtain an $L_{\infty}$ structure on the complex concentrated in degrees $-r+1,\cdots,0$
\begin{equation}\label{complex}
C^{\infty}(M) \overset{d}{\rightarrow}\cdots 
\overset{d}{\rightarrow} \Omega^{r-2}(M)\overset{d}{\rightarrow} \Gamma(E^{r-1})=\chi(M)\oplus \Omega^{r-1}(M),
\end{equation}
for any manifold $M$ and integer $r\ge 2$. When $r=2$ we obtain exactly the Lie 2-algebra given just before Thm. \ref{mor01}.

Let us first recall Getlzer's recent  theorem \cite[Thm. 3]{GetzlerHigherDer}.
Let $(V,\delta,\op\;,\;\cp)$ be a differential graded Lie algebra (DGLA).
  Getlzer  endows the graded\footnote{We take the opposite grading as in \cite{GetzlerHigherDer} 
 so that our differential $\delta$ has degree 1.} vector space
 $V^-:=\oplus_{i<0} V_i$
 with multibrackets satisfying the relations \cite[Def. 1]{GetzlerHigherDer}, which after a degree shift provide 
 $V^-[-1]$ with a $L_{\infty}$-algebra structure in the sense of our Def. \ref{lidef}. Notice that $V^-[-1]$ is concentrated in non-positive degrees: its degree $0$  component is $V_{-1}$, its degree $-1$  component is $V_{-2}$,  and so on.
  The multibrackets are built out of a derived bracket construction using the restriction of the operator $\delta$ to $V_0$, and 
 the Bernoulli numbers appear as coefficients. 

Now let $M$ be a manifold, fix an integer $r\ge 2$, and consider the graded manifold $$T^*[r]T[1]M$$ (see  \cite{Dima}\cite[\S 2]{AlbICM}\cite{ALbFlRio} for background material on graded manifolds). $T^*[r]T[1]M$
is endowed with a canonical Poisson structure of degree $-r$:
there is a  bracket 
$\op\;,\;\cp$ of degree $-r$ on  the graded commutative algebra of functions  $\cC:=C(T^*[r]T[1]M)$ such that
$$\big(\cC \;,\; \cdot \;,\; \op\;,\;\cp\big)$$ is a \emph{Poisson algebra  of degree $r$}
\cite[Def. 1.1]{GPA}. This means that $\op\;,\;\cp$
defines a (degree zero) graded Lie algebra structure on $\cC[r]$, the graded vector space defined by the degree shift $(\cC[r])_i:=\cC_{r+i}$, and that
  $\op a,\cdot \cp$ is a degree $|a|-r$ derivation of the product for any homogeneous element $a\in \cC$.

More concretely, choose coordinates $x_i$  on $M$, inducing fiber coordinates $v_i$ on $T[1]M$, and conjugate coordinates $P_i$ and $p_i$ on the fibers of $T^*[r]T[1]M\to T[1]M$. The degrees of these generators of $\cC$ are $$|x_i|=0,\;\; |v_i|=1,\;\; |P_i|=r,\;\;|p_i|=r-1.$$ Then 
\begin{align*}
\op P_i,x_i \cp=&1=- \op x_i,P_i \cp\\
\op p_i,v_i \cp=&1=-(-1)^{r-1}\op v_i,p_i \cp
\end{align*}
 for all $i$, and all the other brackets between generators vanish. 
Notice that the coordinate $v_i$ corresponds canonically to  $dx_i\in \Omega^1(M)$ and that $p_i$ corresponds canonically to $\pd{x_i}\in \chi(M)$.
Also notice that   $\cC$ is concentrated in non-negative degrees, and that there are canonical identifications 
\begin{equation}\label{ccforms}
\cC_i=\Omega^i(M) \text{ for }0\le i< r-1, \;\;\;\;\;\;\;\;\; \cC_{r-1}=\Omega^{r-1}(M)\oplus \chi(M).
\end{equation}
Indeed
for $i<r-1$
the elements of degree $i$ are sums of expressions of the form $f(x)v_{j_1}\dots v_{j_i}$, while for $i=r-1$ they are sums of expressions $f(x)v_{j_1}\dots v_{j_{r-1}}+g(x)p_j$.
 
 The degree $r+1$ function $\cS:=\sum v_iP_i$, given by the De Rham differential on $M$, satisfies $\op \cS,\cS \cp=0$, hence $\op \cS,\;\cp$ squares to zero. This and the fact that $(\cC[r],\op\;,\;\cp)$ is a graded Lie algebra imply that 
\begin{equation}\label{DGLAme}
\big(\cC[r], \delta:=\op \cS,\;\cp, \op\;,\;\cp\big).
\end{equation}
is a DGLA.
Hence Getlzer's construction can be applied to \eqref{DGLAme}, endowing $(\cC[r])^-[-1]=(\oplus_{0\le i \le  r-1}\cC_i)[r-1]$ (the complex displayed in \eqref{complex}) with an
 $L_{\infty}$-algebra structure.
 
We write out explicitly the multibrackets. The twisted case will be considered in Prop. \ref{ordH} below.
\begin{prop}\label{ord}
Let $M$ be a manifold, $r\ge 2$ an integer.  There exists a Lie $r$-algebra structure on the   complex \eqref{complex} concentrated in degrees $-r+1,\cdots,0$, that is
\begin{equation*} 
C^{\infty}(M) \overset{d}{\rightarrow}\cdots 
\overset{d}{\rightarrow} \Omega^{r-2}(M)\overset{d}{\rightarrow} \Gamma(E^{r-1})=\chi(M)\oplus \Omega^{r-1}(M),
\end{equation*} 
whose only non-vanishing brackets (up to permutations of the entries) are  \begin{itemize}
\item unary bracket: the de Rham differential in negative degrees.

\item  binary bracket: 
\begin{itemize}
\item[ ] for $e_i\in\Gamma(E^{r-1})$  the Courant bracket as in eq. \eqref{dorcou},
$$[e_1,e_2]=\oo e_1,e_2 \cc_{Cou};$$
\item[ ]  for $e=(X,\alpha) \in \Gamma(E^{r-1})$ and $\xi \in \Omega^{\bullet<{r-1}}(M)$,
  $$[e,\xi]=\frac{1}{2} \cL_{X}\xi.$$ \end{itemize}

\item trinary bracket: 
\begin{itemize}
\item[ ]   for   $e_i\in \Gamma(E^{r-1})$,
$$[e_0,e_1,e_2]=-\frac{1}{6}\left( \la \oo e_0,e_1 \cc_{Cou},e_2\ra + c.p.\right);$$   

\item[ ]  for $\xi \in\Omega^{\bullet<{r-1}}(M)$ and $e_i=(X_i,\alpha_i) \in \Gamma(E^{r-1})$,
\begin{align*}
[\xi,e_1,e_2]
=& -\frac{1}{6}\left(
\frac{1}{2}(\iota_{X_1}\cL_{X_2} -  \iota_{X_2}\cL_{X_1}) +\iota_{[X_1,X_2]} \right) \xi.
\end{align*}
\end{itemize}

\item $n$-ary bracket for $n\ge 3$ with $n$ an \emph{odd} integer:  
\begin{itemize}
\item[ ] for $e_i=(X_i,\alpha_i) \in \Gamma(E^{r-1})$, $[e_0,\cdots,e_{n-1}]=\sum_i  [X_0,\dots,\alpha_i,\dots,X_{n-1}]$, with
  \begin{align*}
[\alpha,X_1,\dots,X_{n-1}]=
\frac{(-1)^{\frac{n+1}{2}} 12 B_{n-1} }{(n-1)(n-2)}
 \sum_{1\le i<j\le n-1}(-1)^{i+j+1}\iota_{X_{n-1}}\dots   
  \widehat{\iota_{X_{j}}}\dots \widehat{\iota_{X_{i}}}\dots
\iota_{X_{1}} [\alpha,X_i,X_j];
\end{align*}

\item[ ]for $\xi \in \Omega^{\bullet<{r-1}}(M)$  and $e_i=(X_i,\alpha_i) \in \Gamma(E^{r-1})$,
 \begin{align*}
  [\xi,e_1,\cdots,e_{n-1}]=
 \frac{(-1)^{\frac{n+1}{2}} 12 B_{n-1} }{(n-1)(n-2)}
 \sum_{1\le i<j\le n-1}(-1)^{i+j+1}\iota_{X_{n-1}} 
\dots \widehat{\iota_{X_{j}}}\dots \widehat{\iota_{X_{i}}}\dots
\iota_{X_{_{1}}} [\xi,X_i,X_j].
\end{align*}  
\end{itemize}
\end{itemize}

Here the $B$'s denote the Bernoulli numbers.
\end{prop}

\begin{remark}
 Bering \cite[\S 5.6]{Bering} shows that the vector fields and differential forms on a manifold $M$ are naturally endowed with multibrackets forming an algebraic structure which generalizes $L_{\infty}$-algebras: the quadratic relations satisfied by Bering's multibrackets have Bernoulli numbers as coefficients.  
The multibrackets appearing in Prop. \ref{ord} are similar to Bering's, and
they differ  not only in the coefficients, but also in that the expression for $[\xi,e_1,\cdots,e_{n-1}]$ (for $n\ge 3$) does not appear among Bering's brackets. This is 
a consequence of the fact that Getzler's multibracket are constructed not out of $\delta$, but  out of its restriction to $V_0$.\end{remark} 
 
\begin{remark}
We write more explicitly the trinary bracket of elements $e_i=(X_i,\alpha_i)\in \Gamma(E^{r-1})$: we have
$[e_0,e_1,e_2]=[\alpha_0,X_1,X_2]-[\alpha_1,X_0,X_2]+[\alpha_2,X_0,X_1]$ with 
\begin{align*}
[\alpha_0,X_1,X_2]=
 -\frac{1}{6}\left(
\frac{1}{2}(\iota_{X_1}\cL_{X_2} -  \iota_{X_2}\cL_{X_1})+\iota_{[X_1,X_2]} +\iota_{X_1}\iota_{X_2}d\right)\alpha_0.
\end{align*}
\end{remark}

\begin{proof}
Let $X_1,X_2,\dots \in \chi(M)$ and $\xi_1,\xi_2,\dots$ be differential forms on $M$. In the following we  identify them with elements of $\cC$
as indicated in eq. \eqref{ccforms}, and we adopt the notation introduced in
 the text before Prop.  \ref{ord}.
 The following holds:
 \begin{itemize}
\item[a)] If  $\xi_i\in \Omega^{k_i}(M)$ for $k_1,k_2$ arbitrary, we have  $$\op   X_1+\xi_1  , X_2+\xi_2  \cp =  \iota_{X_1}\xi_2+(-1)^{r-1-k_1} \iota_{X_2}\xi_1.$$
In particular, when  $\xi_1,\xi_2 \in \Omega^{r-1}(M)$, we obtain
  the pairing $\la  \cdot  , \cdot \ra$ as in  eq. \eqref{pairing}.

\item[b)]  For any differential form $\xi_1$,
the identity $$\op\cS,  \xi_1 \cp =d \xi_1$$ is immediate in  coordinates. 

\item[c)] If $\xi_1,\xi_2 \in \Omega^{r-1}(M)$ we have
 $$\op \op\cS, X_1+\xi_1 \cp, X_2+\xi_2  \cp = \oo X_1+\xi_1, X_2+\xi_2 \cc,$$
 the Dorfman bracket as in eq. \eqref{dorf}. 
This holds by the 
following identities, which we write for $\xi_i \in \Omega^{k_i}(M)$ for arbitrary $k_1,k_2$:
$$\op \op\cS, X_1 \cp, X_2  \cp=[X_1,X_2] \text{   and    }\op \op\cS, \xi_1 \cp, \xi_2  \cp=0$$ are checked in coordinates,  and
\begin{align*}
\op \op\cS, X_1  \cp,  \xi_2  \cp&=\op \cS, \op X_1  ,  \xi_2 \cp \cp
+ \op  X_1 , \op \cS, \xi_2 \cp \cp=d(\iota_{X_1} \xi_2)+\iota_{X_1} d\xi_2=\cL_{X_1} \xi_2,\\
\op \op\cS, \xi_1  \cp,  X_2  \cp&=-(-1)^{r-1-k_1}\op   X_2 , \op\cS, \xi_1  \cp \cp=-(-1)^{r-1-k_1}\iota_{X_2} d\xi_1.
\end{align*}

\item[d)] For $n\ge 3$, and letting $a_i$ be either a vector field $X_i$ or a differential form $\xi_i$  of arbitrary degree (not a sum of both),
$$\op \op \dots \op \cS, a_1\cp, \dots \cp,  a_n \cp=0$$ except when 
exactly one of $a_1,a_2,a_3$ is a differential form and all the remaining $a_i$'s are vector fields.
\end{itemize}

Using this it is straighforward to write out the (graded symmetric) multibrackets of \cite[Thm. 3]{GetzlerHigherDer}, which we denote by $(\cdot,\dots,\cdot)$. More precisely, b) gives the unary bracket, c) gives the binary bracket, c) and d) give the trinary bracket.
For the higher brackets ($n\ge 3$ odd) one uses d) and then a) to compute
\begin{align*}
(\alpha,X_1,\dots,X_{n-1})&=
\frac{c_{n-1}}{c_2}\sum_{\sigma \in \Sigma_{n-1}\;,\; \sigma_1 <\sigma_2}
(-1)^{\sigma}\op \op \dots \op (\alpha,X_{\sigma_1},X_{\sigma_2}),
X_{\sigma_3}
\cp, \dots  \cp,  X_{\sigma_{n-1}} \cp\\
&=(-1)^{n-2  \choose 2}(n-3)! \frac{c_{n-1}}{c_2}
\sum_{1\le i<j\le n-1}(-1)^{i+j+1}\iota_{X_{n-1}} 
\dots \widehat{\iota_{X_{j}}}\dots \widehat{\iota_{X_{i}}}\dots
\iota_{X_{_{1}}} (\xi,X_i,X_j),
\end{align*}
where we abbreviate $c_{n-1}:=\frac{(-1)^{n+1\choose 2}}{(n-1)!}B_{n-1}$. The computation for $[\xi,e_1,\cdots,e_{n-1}]$ with $\xi \in \Omega^{\bullet<{r-1}}(M)$   delivers the same expression and uses the fact that $n$ is odd. The coefficient can be simplified:  
$$(-1)^{n-2  \choose 2}(n-3)! \frac{c_{n-1}}{c_2}=\frac{12}{(n-1)(n-2)}B_{n-1}$$
since $n$ is odd and $c_2=\frac{1}{12}$.
 
This gives us the (graded symmetric) multibrackets $(\cdot,\dots,\cdot)$ of  
\cite{GetzlerHigherDer}. As pointed out in \cite{GetzlerHigherDer}, multiplying the   $n$-ary bracket by  
$(-1)^{n-1 \choose 2}$ delivers (graded symmetric) multibrackets 
that satisfy  the Jacobi rules  given just before \cite[Def. 4.2]{GetzlerAnnals}. 

These Jacobi rules coincide with Voronov's \cite[Def. 1]{vor}, and according to  \cite[Rem. 2.1]{vor}, the passage from these (graded symmetric) multibrackets to the (graded skew-symmetric) multibrackets satisfying  our Def. \ref{lidef} is given as follows: multiply 
the multibracket of elements $x_1,\dots,x_n$ by 
\begin{equation}\label{vordec}
(-1)^{\tilde{x}_1(n-1)+\tilde{x}_2(n-2)+\dots+\tilde{x}_{n-1}}
\end{equation}
where $\tilde{x}_i$ denotes the degree of $x_i$ as an element
 of  \eqref{complex}, a complex concentrated in degrees $-r+1,\dots,0$. One easily checks that in all the  cases relevant to us \eqref{vordec} does not introduce any sign.

In conclusion, to pass from the conventions of 
\cite{GetzlerHigherDer} to the conventions of 
our Def. \ref{lidef} we just have to multiply the $n$-ary bracket
$(\cdot,\dots,\cdot)$ 
by $(-1)^{n-1 \choose 2}$,
which for $n=1,2$ equals $1$ and for $n$ odd equals $(-1)^{\frac{n-1}{2}}$.
\end{proof}

Now let  $H\in\Omega^{r+1}_{closed}(M)$ be a closed $r+1$-form. $H$ can be viewed as an element $\cH$ of $\cC_{r+1}$, and  
$\op \cS-\cH,\cS -\cH\cp=-2\op \cS, \cH\cp=-2dH=0$. Hence
\begin{equation}\label{dglaH}
\big(\cC[r], \delta:=\op \cS -\cH,\;\cp, \op\;,\;\cp\big)
\end{equation}
is a DGLA, and again we can apply Getzler's construction. We obtain an $L_{\infty}$-algebra structure that extends the $H$-twisted Courant bracket:
\begin{prop}\label{ordH}
Let $M$ be a manifold, $r\ge 2$ an integer, and $H\in\Omega^{r+1}_{closed}(M)$. There exists a Lie $r$-algebra structure on the   complex \eqref{complex} concentrated in degrees $-r+1,\cdots,0$, whose only non-vanishing brackets (up to permutations of the entries) are those given in Prop. \ref{ord} and additionally 
for $e_i=(X_i,\alpha_i) \in\Gamma(E^{r-1})$:
\begin{itemize}
\item  binary bracket:
$$[e_1,e_2]=\iota_{X_2}\iota_{X_1}H$$
\item $n$-ary bracket for $n\ge 3$ with $n$ an \emph{odd} integer: 
 $$[e_1,\cdots,e_{n}]=(-1)^{\frac{n-1}{2}}\cdot n\cdot B_{n-1}\cdot \iota_{X_n}\dots \iota_{X_1}H.$$
\end{itemize}
\end{prop}
\begin{proof}
It is easy to see (in coordinates, or using that $T[1]M\subset T^*[r]T[1]M$ is Lagrangian) that for any  $n\ge 1$, letting $a_i$ be either a vector field $X_i$ or a differential form $\xi_i$  of arbitrary degree (not a sum of both), one has:
$$\op \op \dots \op \cH, a_1\cp, \dots \cp,  a_n \cp=0$$ except when 
all of the $a_i$'s are vector fields $X_i$'s. In this case one obtains
\begin{equation}\label{Hn}
(-1)^{n \choose 2} \iota_{X_n}\dots \iota_{X_1}H
\end{equation}
using a) in the proof of Prop. \ref{ord}.
Denoting by $(\cdot,\dots,\cdot)$ the (graded symmetric) multibrackets as in \cite{GetzlerHigherDer} from the DGLA \eqref{dglaH}, we see  that   
 $(X_1,\dots,X_n)$   is equal to \eqref{Hn} multiplied by  $-n!\cdot c_{n-1}$.
 In order to pass from the conventions of \cite{GetzlerHigherDer} to those of our Def. \ref{lidef} we multiply by $(-1)^{n-1 \choose 2}$ and obtain the formulae in the statement. \end{proof}

For any $B\in \Omega^r(M)$, the  gauge
transformation of $E^{r-1}$ given
by $e^{-B} \colon X+\alpha \mapsto X+\alpha-\iota_XB$ maps the $H$-twisted Courant bracket to the $(H+dB)$-twisted Courant bracket. Defining properly the notion of higher Courant algebroid -- of which the $E^r$'s should be the main examples -- and
extending to this general setting   Prop. \ref{ord}, will presumably  imply that
  the $L_{\infty}$-algebras defined  by cohomologous differential forms are isomorphic. We show this directly:
\begin{prop}
 Let $M$ be a manifold, $r\ge 2$ an integer, and $H\in\Omega^{r+1}_{closed}(M)$. For any $B\in \Omega^r(M)$, there is a strict isomorphism 
 \begin{center}
$($the Lie-$r$ algebra defined  by $H)\to($the Lie-$r$ algebra defined  by $H+dB)$
   \end{center}
between the Lie-$r$ algebra structures defined as in 
 Prop. \ref{ordH}
on the complex \eqref{complex}. Explicitly, the isomorphism is given by $e^{-B}$ on $\Gamma(E^{r-1})$ and is the identity elsewhere. 
\end{prop}
\begin{proof}
View $B$ as an element $\cB\in\cC_r$. As $\op\cB,\;\cp$ is a degree zero derivation of the graded Lie algebra
$(\cC[r],  \op\;,\;\cp)$ and is nilpotent, it follows that the exponential $\Phi:=e^{\op\cB,\;\cp}$ is an automorphism. Therefore it is an isomorphism   of  DGLAs $$\Phi \colon \big(\cC[r], \delta:=\op \cS -\cH,\;\cp, \op\;,\;\cp\big) \to \big(\cC[r], \Phi\delta \Phi^{-1}, \op\;,\;\cp\big).$$
From the formulas for the multibrackets in Getzler's \cite[Thm. 3]{GetzlerHigherDer} it is then clear that $\Phi|_{(\oplus_{0\le i \le  r-1}\cC_i)[r-1]}$ is a strict isomorphism between the  $L_{\infty}$-algebras induced by these two DGLAs.

The differential $\Phi\delta \Phi^{-1}$ on $\cC$ is not equal to  $\op \cS - (\cH+ \op \cS,\cB \cp),\;\cp$, which is the differential associated to     $H+dB\in\Omega^{r+1}_{closed}(M)$ as in \eqref{dglaH}. However
on $\oplus_{0\le i \le  r-1}\cC_i$ the two differentials do agree.
(This follows from the fact that on $\oplus_{0\le i \le  r-1}\cC_i$ we have $\Phi(y) =y+\op\cB,y\cp$).    
This assures that the $L_{\infty}$-algebras induced by the two differentials agree.
\end{proof}

\section{Open questions: the relation between the $L_{\infty}$-algebras of \S\ref{obs} and \S\ref{liM}-\ref{ez}}\label{per}

In this section we speculate about the relations among the $L_{\infty}$-algebras that appeared in \S \ref{obs}--\S\ref{ez} and their higher analogues, and relate them to prequantization.

Let $M$ be a manifold. Given an   integer $n\ge0$ and $H\in \Omega^{n+2}_{closed}(M)$, we use the notation $E^n_H$ to denote the vector bundle
$E^n=TM\oplus \wedge^n T^*M$ with the $H$-twisted Dorfman bracket $[\cdot,\cdot]_{H}$.
In particular, $E^n_0$ denotes $TM\oplus \wedge^n T^*M$ with the untwisted Dorfman bracket \eqref{dorf}.

\subsection{Relations between $L_{\infty}$-algebras}\label{relations}
To any $n\ge 0$ and $H\in \Omega^{n+2}_{closed}(M)$, we associated in Prop. \ref{ordH}  
    a Lie $n+1$-algebra
$\cS^{E^n_H}$. We ask:
\begin{itemize} 
\item[]
 Is there a natural $L_{\infty}$-morphism $D$ from 
$ \cS^{E^n_H}$ to $\cS^{E_0^{n+1}}$?
\end{itemize} 
When $n=0$ the answer is affirmative by Thm. \ref{mor01}.

Let $p\ge1$ and  
$L\subset E^{p}_0$ an involutive isotropic subbundle. Denote by $\cO^{L\subset E^p_0}$ the Lie $p$-algebra associated  in Thm. \ref{Liep}.
  Since $L$ is an involutive subbundle of  $E_0^p$ it is natural to ask:
\begin{itemize} 
\item[]  What is the relation between $\cO^{L\subset E^p_0}$ 
and $\cS^{E_0^{p}}$?
\end{itemize}
When $L$ is equal to $graph({H})$ for a  $p$-plectic form ${H}$, we expect the relation 
to be  given by  an $L_{\infty}$-morphism
 \begin{equation*} 
P\colon \cO^{graph({H})\subset E^p_0} \rightsquigarrow \cS^{E_H^{p-1}}
\end{equation*}
with the property that the unary map of  the
 $L_{\infty}$-morphism $D \circ P$, restricted to the degree zero component,   coincide with 
\begin{equation}\label{sortapreqmap}
{\Omega}^{p-1}_{ham}(M, graph({H})) \to \Gamma(E^p_0), \;\;\;\;\;\alpha \mapsto X_{\alpha}-d\alpha. 
\end{equation}
We summarize the situation in this diagram:\[\xymatrix{
\cS^{E^{p-1}_H}\ar@{~>}[r]^D &    \cS^{E^{p}_0}\\
\ar@{~>}[u]^P\ar@{~>}[ur]_{D \circ P}\cO^{graph({H})\subset E^p_0}  & 
}
\]

\begin{remark}\label{expre}
In the case $p=1$ (so ${H}$ is a symplectic form) the embedding $P$ exists and is given as follows. We have two honest Lie algebras 
 $$
\cO^{graph({H})\subset E^1_0}=(C^{\infty}(M),\{\cdot,\cdot \}),\;\;\;\;\;\; 
\cS^{E_H^0}=(\Gamma(TM\oplus \RR),[\cdot,\cdot]_{H})
$$
where $\{\cdot,\cdot\}$ is the usual Poisson bracket defined by $H$.
The map  
$$
P \colon C^{\infty}(M) \to \Gamma(TM\oplus \RR),\;\;\;\;
 f \mapsto (X_f,-f)
$$
 is a Lie algebra morphism. Lie $2$-algebra morphism $D$ is given by  Thm. \ref{mor01}. One computes that the composition consists only of a unary map, given by the Lie algebra morphism \eqref{sortapreqmap}.
\end{remark}

\begin{remark} We interpret $P$ as a  prequantization map. Indeed for $p=1$ and integral form $H$,  the Lie algebra $\cS^{E_H^{0}}$ can be identified with the space of $S^1$-invariant vector fields on a circle bundle over $M$ \cite[\S 3.8]{Gu}. The composition of $P$  with the  action of vector fields on the $S^1$-equivariant complex valued functions is then a   faithful representation of the Lie algebra $\cO^{graph({H})\subset E^1_0}=C^{\infty}(M)$, that is, a prequantization representation.  For $p=2$ the morphism $P$ is described by  Rogers in  \cite[Thm. 5.2]{RogersCou} and   \cite[Thm. 7.1]{RogersPre}, to which we refer for the interpretation as a prequantization map. 
\end{remark}

\subsection{The twisted case} We pose three questions about higher analogues of twisted Dirac structures.
Let $H$ be a closed $p+1$-form for $p\ge 2$. Let $L'\subset E^{p-1}_H$
be an isotropic subbundle, involutive w.r.t. the ${H}$-twisted Dorfman bracket.

\begin{itemize}
\item[] 
Can one associate to $L'$ an $L_{\infty}$-algebra of  observables  
$\cO^{L'\subset E^{p-1}_H}$?
\end{itemize}
To the author's knowledge, this is not known even in the simplest case, i.e., when $p=2$ and  $L'$ is the graph of an  $H$-twisted Poisson structure \cite{SW}. In that case one  defines in the usual manner a skew-symmetric  bracket $\{\cdot,\cdot\}$ on $C^{\infty}(M)$. It does not satisfy the Jacobi identity but rather \cite[eq. (4)]{SW}
$\{\{f,g\},h\} +c.p.=-{H}(X_f,X_g,X_h)$, hence it is natural to wonder if 
one can extend
this bracket to an $L_{\infty}$-structure.

\begin{itemize}
\item []
Is there a natural $L_{\infty}$-morphism $D'$ from  $\cO^{L'\subset E^{p-1}_H}$ to $\cO^{graph(H)\subset E^{p}_0}$?
\end{itemize}
This question is motivated by the fact that $L'$ plays the role of a primitive of $H$. In the simple case that $L'$ is the graph of a symplectic form the answer is affirmative, by
 the morphism from $(C^{\infty}(M),\{\cdot,\cdot\})$ to $C^{\infty}(M) \overset{d}{\rightarrow}\Omega^1_{closed}(M)$ (a complex with no higher brackets) with  vanishing unary map and   binary map $\phi_2(f,g)=\{f,g\}$.  

\begin{itemize}
\item  []
 Is there an $L_{\infty}$-morphism from
$ \cO^{L'\subset E^{p-1}_H}$ to $\cS^{E^{p-1}_H}$, assuming that $L'$ is the graph of a non-degenerate differential form?
\end{itemize}
Such a morphism would be interesting because it could be interpreted as a
 weaker (because not injective) version of a prequantization map for $(M,L')$. \\

We summarize the  discussion of this whole section in the following diagram, in which for the sake of concreteness and simplicity we take ${H}\in \Omega^3_{closed}(M)$ to be  a 2-plectic form and $L'\subset TM\oplus T^*M$ to be a ${H}$-twisted Dirac structure. The arrows denote $L_{\infty}$-morphisms.

\[\xymatrix{
& \cS^{E^{1}_H} \ar@{~>}[r]^D &     \cS^{E^{2}_0}\\
\cO^{L'\subset E^{1}_H} \ar@{~>}[r]^{D'}\ar@{~>}[ur] &\ar@{~>}[u]^P\ar@{~>}[ur]_{D \circ P} \cO^{graph({H})\subset E_0^2}  & 
}
\]


We conclude presenting an interesting  example in which the geometric set-up  described above applies.
\begin{ex} 
Let $G$ be  a Lie group whose   Lie algebra $\g$  is endowed with a non-degenerate bi-invariant quadratic form $(\cdot,\cdot)_{\g}$. There is  a well-defined closed Cartan 3-form $H$, which on $\g=T_eG$ is given by $H(u,v,w)=\frac{1}{2} 
(u,[v,w])_{\g}$ \cite[\S 2.3]{MaHenr}. 
There is also a canonical  $H$-twisted Dirac structure  $L'\subset TG\oplus T^*G$: it is    given by  $L'=\{(v_r-v_l)+\frac{1}{2}(v_r+v_l)^*:v\in \g \}$ where $v_r,v_l$ denote the right and left translations of $v\in \g$ and the quadratic form is used to identify a tangent vector $X\in TG$ with a covector $X^*\in T^*G$
 \cite{SW}\cite[Ex. 3.4]{MaHenr}. \end{ex}

\appendix

\section{The proof of Proposition \ref{linalg}}\label{applag}

In this appendix we present the proof of Prop. \ref{linalg}. We start giving an alternative characterization of Lagrangian subspaces.

\begin{lemma}\label{easychar}
Let $T$ be a vector space and $p\ge 1$.
For all subspaces $L\subset T\oplus \wedge^p T^*$, denoting $S:=pr_T L$, the following holds:
$$L \text{ is Lagrangian }\Leftrightarrow 
\begin{cases}
 L \text{ is isotropic}\\
 L\cap \wedge^p T^*=\wedge^p S^{\circ}\\
 dim(S)\le (dim(T)-p) \text{ or }S=T. \end{cases}$$
\end{lemma}
\begin{proof}
``$\Rightarrow$:'' 
Assume first that $L$ is Lagrangian. It is straightforward to check that for any subspace $F\subset T\oplus \wedge^p T^*$ we have
\begin{equation}\label{ker}
F^{\perp}\cap \wedge^pT^*=\wedge^p (pr_T(F))^{\circ}.
\end{equation}
We apply this to 
$F=L=L^{\perp}$ and derive $L\cap \wedge^p T^*=\wedge^p S^{\circ}$. 

Hence we are left with  showing that $S$ satisfies $dim(S)\le dim(T)-p \text{ or }S=T$.
We argue by contradiction: we 
assume that
$\wedge^p S^{\circ}= \{0\}$ and $S$ is strictly included in $T$, and deduce from this  that $pr_T(L^{\perp})\not\subset S$, which contradicts
$L= L^{\perp}$. 
Let $\{X_j\}_{j\le dim(T)}$ be a basis of $T$ whose first $dim(S)$ elements form a basis of $S$. Let $Y$ be a basis element \emph{not} lying in $S$ (it exists since $S\neq T$).  It is enough to prove 
 the following \emph{claim}:
\begin{equation*}Y+\beta \in L^{\perp} \text{ 
where }\;\;\;
\beta=-\sum_{j=1}^{dim(S)}\left(X^*_j\wedge(\sum_{q=0}^p \frac{1}{q+1}\iota_Y \alpha_j^q)\right),
\end{equation*}
because it implies that  $Y\in pr_T(L^{\perp})$.
Here $\{X^*_j\}_{j\le dim(T)}$ denotes the  basis of $T^*$ dual to $\{X_j\}_{j\le dim(T)}$,
and, for all $j\le dim(S)$, $\alpha_j\in \wedge^p T^*$ is such that $X_j+\alpha_j \in L$. Further we adopt the following notation: 
for any $\alpha \in \wedge^p T^*$, $\alpha^q$  denotes the component of $\alpha$,  written in the basis of $\wedge^p T^*$ induced by  $\{X^*_j\}_{j\le dim(T)}$,  for which  the number of $X^*_j$'s with $j\le dim(S)$ is exactly $q$.

To prove the claim fix $j_0\le dim(S)$. We have
\begin{align*}
\iota_{X_{j_0}}\beta 
&=- \sum_{q=0}^p \frac{1}{q+1}\iota_Y \alpha_{j_0}^q
+\sum_{j=1}^{dim(S)}X^*_j\wedge(\sum_{q=0}^p \frac{1}{q+1}\iota_{X_{j_0}}\iota_Y \alpha_j^q)\\
&=- \sum_{q=0}^p \frac{1}{q+1}\iota_Y \alpha_{j_0}^q
- \iota_Y\sum_{j=1}^{dim(S)}X^*_j\wedge (\sum_{q=0}^p \frac{1}{q+1}  \iota_{X_{j}}\alpha_{j_0}^q)\\
&=- \sum_{q=0}^p \frac{1}{q+1}\iota_Y \alpha_{j_0}^q
-  \sum_{q=0}^p \frac{q}{q+1} \iota_Y \alpha_{j_0}^q\\
&=-  \iota_Y \alpha_{j_0},
\end{align*}
where in the second equality we used
$\iota_{X_{j_0}}\alpha_j^q=-\iota_{X_{j}}\alpha_{j_0}^q$
and in the third
$\sum_{j=1}^{dim(S)}X^*_j\wedge(   \iota_{X_{j}} \alpha_{j_0}^q)=q\alpha_{j_0}^q$.
Hence $\la Y+\beta, X_j+\alpha_j \ra =0$ for all $j\le dim(S)$. 
Since 
$L\cap \wedge^p T^*=\wedge^p S^{\circ}= \{0\}$, we have $L=span 
 \{X_j+\alpha_j\}_{j\le dim(S)}$, and we conclude that 
  $Y+\beta \in L^{\perp}$, proving the claim.

``$\Leftarrow$:'' 
We need to show that $L$ is Lagrangian, i.e. $L=L^{\perp}$.
We claim that  $pr_T(L^{\perp})=S$.  If $S=T$ this is clear, so we prove the claim
 in the case $dim(S)\le dim(T)-p$, for which we have $\wedge^p S^{\circ}\neq \{0\}$. Since  $\wedge^p S^{\circ}\subset L$, this implies that $pr_T(L^{\perp})\subset S$. By the isotropicity of $L$ we therefore have
 $pr_T(L^{\perp})=S$, as claimed.
  
Hence if $X+\beta \in L^{\perp}$ there exists $\alpha\in \wedge^{p}T^*$ such that
$X+\alpha \in L\subset L^{\perp}$. So $\beta-\alpha\in L^{\perp}\cap \wedge^{p}T^*=\wedge^p S^{\circ}\subset L$, where the equality holds by eq. \eqref{ker}. Therefore
$X+\beta=(X+\alpha)+(\beta-\alpha)$ is the sum of two elements of $L$, showing $L^{\perp}\subset L$.
\end{proof}

\begin{lemma}\label{ext} Let $S\subset T$ a subspace and $p\ge 1$.
Let $\Omega \in \wedge^2 S^*\otimes \wedge^{p-1}T^*$. Then $\Omega$ admits an extension to   $S^*\otimes \wedge^{p}T^*$ if{f} it admits an
extension to   $ \wedge^{p+1}T^*$.
\end{lemma}
\begin{proof}
If there exists  $\alpha\in \wedge^{p+1}T^*$ with $\alpha|_{S\otimes S\otimes\bigotimes^{p-1}T}=\Omega$, the clearly $\alpha|_{S\otimes \bigotimes^{p}T}$ is an element of $S^*\otimes \wedge^{p}T^*$ with the required property.

Conversely, let $\beta' \in S^*\otimes \wedge^{p}T^*$ be an extension of $\Omega$. We choose a complement $C$ to $S$ in $T$, and by the identification $S^*\cong C^{\circ}$ from $\beta'$ we obtain an element $\beta\in  T^*\otimes \wedge^{p}T^*$. The skew-symmetrization  $\bar{\beta} \in  \wedge^{p+1}T^*$ of $\beta$ is given as follows:   $$\bar{\beta}(x_0,\dots,x_{p})=\frac{1}{p+1}\sum_{j=0}^p
(-1)^j \beta(x_j,x_0,\dots,\hat{x}_{j},\dots,x_{p})$$ for all $x_i\in T$.
In general $\bar{\beta}$ does not restrict to $\Omega$, but a weighted sum of its component does, as we now show.
We have $\bar{\beta}=\sum_{q=0}^{p+1} \bar{\beta}^q$. Here,
for any basis $\{X_j\}_{j\le dim(T)}$ of $T$ whose first $dim(S)$ elements span $S$ and whose remaining elements span $C$, taking $\{X^*_j\}_{j\le dim(T)}$ to be the dual basis of $T^*$, 
we denote by  $ \bar{\beta}^q$  the component of $\bar{\beta}$ for which,  in the basis of $\wedge^p T^*$ induced by  $\{X^*_j\}_{j\le dim(T)}$,    the number of $X^*_j$'s with $j\le dim(S)$ is exactly $q$. We have $ \bar{\beta}^0=0$, since $\beta$ is an extension of $\beta'$. For $q=1,\dots, p+1$, vectors $x_0,\dots,x_{q-1}\in S$ and $x_q,\dots,x_p\in C$ we have\footnote{This of course does not  imply that $\beta$ is totally skew, as the element
 $x_0$ of $S$ is plugged in the first slot of $\beta$.}  
$$ \bar{\beta}^q (x_0,\dots,x_p)=\bar{\beta} (x_0,\dots,x_p)=\frac{q}{p+1}\beta(x_0,\dots,x_p).$$ Therefore $\sum_{q=1}^{p+1} \frac{p+1}{q} \bar{\beta}^q$ is an element of  $\wedge^{p+1}T^*$ whose restriction to $S\otimes \bigotimes^{p}T$ agrees with $\beta'$, and in particular its restriction to $S\otimes S \otimes \bigotimes^{p-1}T$ agrees with $\Omega$.
\end{proof}

\begin{proof}[Proof of Prop. \ref{linalg}]
We make use of the characterization of Lagrangian subspaces given in Lemma \ref{easychar}.

 We first show that the correspondence ``$L \mapsto (S,\Omega)$'' is well-defined. 
Let $L$ be a Lagrangian subspace.
The dimension restriction on $S$ follows from Lemma \ref{easychar}. 
Since $L\cap \wedge^p T^*=\wedge^p S^{\circ}$, for any $X\in S$,  the   definition of  $\iota_X\Omega$ in Prop. \ref{linalg} is independent of 
 the choice of $\alpha$ with $X+\alpha \in L$, and determines a unique $\Omega \in 
 \otimes^2 S^*\otimes \wedge^{p-1}T^*$. Clearly $\Omega$ is skew in the first two components: if $X+\alpha,Y+\beta \in L$ then the isotropicity of $L$ implies
 $\iota_Y\iota_X \Omega=\iota_Y\alpha=-\iota_X\beta=
-\iota_X\iota_Y \Omega$. By construction, $\Omega$ is the restriction of an element of $S^*\otimes \wedge^{p}T^*$, hence by Lemma \ref{ext} it is the restriction of an element of $\wedge^{p+1}T^* $
 
Next, we show that the correspondence ``$(S,\Omega)\mapsto L$'' is well-defined.
Let $(S,\Omega)$ a pair as in the statement of Prop. \ref{linalg}. This pair  maps to  a subspace $L$ which is isotropic, due to the skew-symmetry of $\Omega$ in its first 2 components. By inspection we have $ L\cap \wedge^p T^*=\wedge^p S^{\circ}$, and further $S$ agrees with $pr_T(L)$ because $\Omega$ is the restriction of an element of $S^*\otimes \wedge^{p}T^*$. Hence  $L$ is Lagrangian by Lemma \ref{easychar}.

The maps  ``$L \mapsto (S,\Omega)$'' and ``$(S,\Omega)\mapsto L$''
are   inverses of each other. 
\end{proof}
 

\bibliographystyle{habbrv}
\bibliography{HigherDirac.bib}

\begin{thebibliography}{10}

\bibitem{Abad:2009zr}
C.~A. Abad and M.~Crainic.
\newblock {The Weil algebra and the Van Est isomorphism}.
\newblock 01 2009, Arxiv:\texttt{0901.0322}.

\bibitem{BHR}
J.~C. Baez, A.~E. Hoffnung, and C.~L. Rogers.
\newblock Categorified symplectic geometry and the classical string.
\newblock \emph{Commun.Math.Phys.293}:701-725,2010.

\bibitem{Bering}
K.~Bering.
\newblock Non-commutative {B}atalin-{V}ilkovisky algebras, homotopy {L}ie
  algebras and the {C}ourant bracket.
\newblock {\em Comm. Math. Phys.}, 274(2):297--341, 2007.

\bibitem{MaHenr}
H.~Bursztyn and M.~Crainic.
\newblock Dirac structures, momentum maps, and quasi-{P}oisson manifolds.
\newblock In {\em The breadth of symplectic and {P}oisson geometry}, volume 232
  of {\em Progr. Math.}, pages 1--40. Birkh{\"a}user Boston, Boston, MA, 2005.

\bibitem{BCWZ}
H.~Bursztyn, M.~Crainic, A.~Weinstein, and C.~Zhu.
\newblock Integration of twisted {D}irac brackets.
\newblock {\em Duke Math. J.}, 123(3):549--607, 2004.

\bibitem{BR}
H.~Bursztyn and O.~Radko.
\newblock Gauge equivalence of {D}irac structures and symplectic groupoids.
\newblock {\em Ann. Inst. Fourier (Grenoble)}, 53(1):309--337, 2003.

\bibitem{CW}
A.~Cannas~da Silva and A.~Weinstein.
\newblock {\em Geometric models for noncommutative algebras}, volume~10 of {\em
  Berkeley Mathematics Lecture Notes}.
\newblock American Mathematical Society, Providence, RI, 1999.

\bibitem{CIDL}
F.~Cantrijn, A.~Ibort, and M.~de~Le{{\'o}}n.
\newblock On the geometry of multisymplectic manifolds.
\newblock {\em J. Austral. Math. Soc. Ser. A}, 66(3):303--330, 1999.

\bibitem{Hammulti}
F.~Cantrijn, L.~A. Ibort, and M.~de~Le{\'o}n.
\newblock Hamiltonian structures on multisymplectic manifolds.
\newblock {\em Rend. Sem. Mat. Univ. Politec. Torino}, 54(3):225--236, 1996.
\newblock Geometrical structures for physical theories, I (Vietri, 1996).

\bibitem{AlbICM}
A.~S. Cattaneo.
\newblock From topological field theory to deformation quantization and
  reduction.
\newblock {Proceedings of ICM 2006, Vol. III, 339-365.} {\tt
  http://www.math.uzh.ch/fileadmin/math/preprints/icm.pdf}.

\bibitem{GPA}
A.~S. Cattaneo, D.~Fiorenza, and R.~Longoni.
\newblock {Graded Poisson Algebras}.
\newblock Encyclopedia of Mathematical Physics, eds. J.-P. Fran{\c c}oise, G.L.
  Naber and Tsou S.T. , vol. 2, p. 560-567 (Oxford: Elsevier, 2006).
  \texttt{http://www.math.uzh.ch/fileadmin/math/preprints/15-05.pdf}.

\bibitem{ALbFlRio}
A.~S. Cattaneo and F.~Sch\"atz.
\newblock Introduction to supergeometry.
\newblock {To appear in the ``Poisson 2010'' proceedings}.

\bibitem{Cou}
T.~J. Courant.
\newblock Dirac manifolds.
\newblock {\em Trans. Amer. Math. Soc.}, 319(2):631--661, 1990.

\bibitem{MariusPre}
M.~Crainic.
\newblock Prequantization and {L}ie brackets.
\newblock {\em J. Symplectic Geom.}, 2(4):579--602, 2004.

\bibitem{MarlkDoubek}
M.~Doubek, M.~Markl, and P.~Zima.
\newblock Deformation theory (lecture notes).
\newblock {\em Archivum mathematicum 43(5), 2007, 333-371}, 05 2007,
  ArXiv:\texttt{0705.3719}.

\bibitem{ekzab}
J.~Ekstrand and M.~Zabzine.
\newblock Courant-like brackets and loop spaces.
\newblock ArXiv:\texttt{0903.3215}.

\bibitem{GetzlerAnnals}
E.~Getzler.
\newblock Lie theory for nilpotent {$L_\infty$}-algebras.
\newblock {\em Ann. of Math. (2)}, 170(1):271--301, 2009.

\bibitem{GetzlerHigherDer}
E.~Getzler.
\newblock Higher derived brackets, 10 2010, ArXiv:\texttt{1010.5859v1}.

\bibitem{Gu}
M.~Gualtieri.
\newblock {Generalized complex geometry}, ArXiv:\texttt{math.DG/0401221}.

\bibitem{Gu2}
M.~Gualtieri.
\newblock Generalized complex geometry, ArXiv:\texttt{math/0703298}.

\bibitem{Hagi}
Y.~Hagiwara.
\newblock Nambu-{D}irac manifolds.
\newblock {\em J. Phys. A}, 35(5):1263--1281, 2002.

\bibitem{Hi}
N.~Hitchin.
\newblock Generalized {C}alabi-{Y}au manifolds.
\newblock {\em Q. J. Math.}, 54(3):281--308, 2003.

\bibitem{LadaMarkl}
T.~Lada and M.~Markl.
\newblock Strongly homotopy {L}ie algebras.
\newblock {\em Comm. Algebra}, 23(6):2147--2161, 1995.

\bibitem{LadaStasheff}
T.~Lada and J.~Stasheff.
\newblock Introduction to {SH} {L}ie algebras for physicists.
\newblock {\em Internat. J. Theoret. Phys.}, 32(7):1087--1103, 1993.

\bibitem{LWX}
Z.-J. Liu, A.~Weinstein, and P.~Xu.
\newblock Manin triples for {L}ie bialgebroids.
\newblock {\em J. Differential Geom.}, 45(3):547--574, 1997.

\bibitem{Moser}
J.~Moser.
\newblock On the volume elements on a manifold.
\newblock {\em Trans. Amer. Math. Soc.}, 120:286--294, 1965.

\bibitem{RogerVol}
C.~Roger.
\newblock Unimodular vector fields and deformation quantization.
\newblock In {\em Deformation quantization ({S}trasbourg, 2001)}, volume~1 of
  {\em IRMA Lect. Math. Theor. Phys.}, pages 135--148. de Gruyter, Berlin,
  2002.

\bibitem{RogersCou}
C.~L. Rogers.
\newblock Courant algebroids from categorified symplectic geometry, 2009,
  ArXiv:\texttt{1001.0040}.

\bibitem{RogersPre}
C.~L. Rogers.
\newblock {2-plectic geometry, Courant algebroids, and categorified
  prequantization}.
\newblock 09 2010, ArXiv\texttt{1009.2975v1}.

\bibitem{RogersL}
C.~L. Rogers.
\newblock {$L_{\infty}$ algebras from multisymplectic geometry}.
\newblock 05 2010, ArXiv:\texttt{1005.2230v1}.

\bibitem{Dima}
D.~Roytenberg.
\newblock On the structure of graded symplectic supermanifolds and {C}ourant
  algebroids.
\newblock In {\em Quantization, Poisson brackets and beyond (Manchester,
  2001)}, volume 315 of {\em Contemp. Math.}, pages 169--185. Amer. Math. Soc.,
  Providence, RI, 2002.

\bibitem{rw}
D.~Roytenberg and A.~Weinstein.
\newblock Courant algebroids and strongly homotopy {L}ie algebras.
\newblock {\em Lett. Math. Phys.}, 46(1):81--93, 1998.

\bibitem{MultiDirac}
J.~Vankerschaver, H.~Yoshimura, and J.~E. Marsden.
\newblock {Multi-Dirac Structures and Hamilton-Pontryagin Principles for
  Lagrange-Dirac Field Theories}.
\newblock 08 2010, Arxiv:\texttt{1008.0252v2}.

\bibitem{vor}
T.~Voronov.
\newblock Higher derived brackets and homotopy algebras.
\newblock {\em J. Pure Appl. Algebra}, 202(1-3):133--153, 2005.

\bibitem{SW}
P.~\v{S}evera and A.~Weinstein.
\newblock Poisson geometry with a 3-form background.
\newblock {\em Progr. Theoret. Phys. Suppl.}, (144):145--154, 2001.

\end{thebibliography}
\end{document}